\documentclass[a4paper]{amsart}
\usepackage{amsmath,amsthm,amssymb}

\newcommand{\deom}{\Delta_\Omega}
\newcommand{\Wom}{W_\Omega}
\newcommand{\Omtil}{\tilde{\Omega}}
\newcommand{\M}{\operatorname{M}}
\newcommand{\ind}{\operatorname{Index}}
\newcommand{\Trace}{\operatorname{Tr}}
\newcommand{\omtil}{\tilde{\omega}}
\newcommand{\Ntil}{\tilde{N}}
\newcommand{\Etil}{\tilde{E}}
\newcommand{\real}{^{\text{\rm\tiny real}}}
\newcommand{\Ah}{\hat{A}}
\newcommand{\uni}{_{\text{\rm u}}}

\newcommand{\vareps}{\varepsilon}
\newcommand{\zetatil}{\tilde{\zeta}}
\newcommand{\etatil}{\tilde{\eta}}
\newcommand{\gatil}{\tilde{\gamma}}
\newcommand{\mutil}{\tilde{\mu}}
\newcommand{\Hbol}{{\overset{\circ}{H}}}
\newcommand{\Honebol}{{\overset{\circ}{H}_1}}
\newcommand{\Htwobol}{{\overset{\circ}{H}_2}}

\newcommand{\F}{\mathbb{F}}
\newcommand{\free}[1]{\Gamma(#1)^{\prime\prime}}
\newcommand{\Cl}{\operatorname{Cl}(H_\R)}
\newcommand{\sitil}{\tilde{\sigma}}
\newcommand{\nabte}{\nabla_\te}
\newcommand{\Sh}{\hat{S}}
\newcommand{\B}{\operatorname{B}}
\newcommand{\site}{\sigma^\theta}
\newcommand{\nabtil}{\tilde{\nabla}}
\newcommand{\el}{\ell}
\newcommand{\tauh}{\hat{\tau}}
\newcommand{\Ad}{\operatorname{Ad}}
\newcommand{\SL}{\operatorname{SL}}

\newcommand{\cA}{\mathcal{A}}

\newcommand{\kruisje}[1]{\, \mbox{$_{#1}$}\hspace{-.2ex}\mbox{$\ltimes$} \,}

\newcommand{\Mh}{\hat{M}}
\newcommand{\Si}{\Sigma}

\newcommand{\dpr}{^{\prime\prime}}

\newcommand{\te}{\theta}
\newcommand{\tetil}{\tilde{\theta}}

\newcommand{\Tr}{{\operatorname{Tr}}}

\newcommand{\cWh}{{\hat{\cW}}}

\newcommand{\Mtil}{\tilde{M}}

\newcommand{\cR}{\mathcal{R}}
\newcommand{\Rh}{\hat{R}}
\newcommand{\R}{\mathbb{R}}
\newcommand{\Z}{\mathbb{Z}}
\newcommand{\C}{\mathbb{C}}
\newcommand{\N}{\mathbb{N}}

\newcommand{\Ga}{\Gamma}
\newcommand{\deh}{\hat{\Delta}}

\newcommand{\vfih}{\hat{\vfi}}

\newcommand{\cU}{\mathcal{U}}

\newcommand{\cJ}{\mathcal{J}}
\newcommand{\nab}{\nabla}
\newcommand{\Jh}{\hat{J}}
\newcommand{\nabh}{\hat{\nab}}
\newcommand{\cZ}{\mathcal{Z}}

\newcommand{\cL}{\mathcal{L}}
\newcommand{\cF}{\mathcal{F}}

\newcommand{\cM}{\mathcal{M}}

\newcommand{\ot}{\otimes}
\newcommand{\sla}{\lambda}
\newcommand{\la}{\Lambda}
\newcommand{\Om}{\Omega}
\newcommand{\om}{\omega}
\newcommand{\io}{\iota}
\newcommand{\vfi}{\varphi}
\newcommand{\eps}{\varepsilon}

\newcommand{\al}{\alpha}
\newcommand{\be}{\beta}
\newcommand{\ga}{\gamma}
\newcommand{\sde}{\delta}
\newcommand{\de}{\Delta}

\newcommand{\si}{\sigma}
\newcommand{\Mfi}{\mathcal{M}_{\vfi}}
\newcommand{\Nfi}{\mathcal{N}_{\vfi}}
\newcommand{\Mpsi}{\mathcal{M}_{\psi}}

\newcommand{\cV}{\mathcal{V}}
\newcommand{\cW}{\mathcal{W}}
\newcommand{\cK}{\mathcal{K}}

\newcommand{\dehop}{\deh \hspace{-.3ex}\raisebox{0.9ex}[0pt][0pt]{\scriptsize\fontshape{n}\selectfont op}}
\newcommand{\recht}{\rightarrow}
\newcommand{\Gammah}{\hat{\Gamma}}
\newcommand{\Ahu}{\hat{A}\uni}
\newcommand{\dehu}{\hat{\Delta}\uni}
\newcommand{\epsh}{\hat{\varepsilon}}
\newcommand{\Rhu}{\hat{R}\uni}
\newcommand{\Wh}{\hat{W}}
\newcommand{\Au}{A\uni}
\newcommand{\deu}{\Delta\uni}
\newcommand{\Shu}{\hat{S}\uni}
\newcommand{\tauhu}{\hat{\tau}^{\text{\rm u}}}

\numberwithin{equation}{section}

{\theoremstyle{definition}\newtheorem{definition}{Definition}[section]

\newtheorem{remark}[definition]{Remark}
}
\newtheorem{proposition}[definition]{Proposition}
\newtheorem{lemma}[definition]{Lemma}
\newtheorem{theorem}[definition]{Theorem}
\newtheorem{corollary}[definition]{Corollary}

\renewcommand{\theenumi}{\alph{enumi}}

\newcommand{\adres}{{\footnotesize   Institut de Math{\'e}matiques de
  Jussieu, Alg{\`e}bres d'Op{\'e}rateurs et
Repr{\'e}sentations \\ 175, rue du Chevaleret; F--75013 Paris (France) \\
  e-mail: vaes@math.jussieu.fr}}

\begin{document}
\title{Strictly outer actions of groups and quantum groups}
\author[Stefaan Vaes]{Stefaan Vaes \vspace{0.5cm} \\ \protect\adres}
\address{Institut de Math{\'e}matiques de Jussieu \\ Alg{\`e}bres d'Op{\'e}rateurs et
Repr{\'e}sentations \newline 175, rue du Chevaleret \\ F--75013 Paris (France)}
\email{vaes@math.jussieu.fr}

\begin{abstract}\noindent
An action of a locally compact group or quantum group on a factor is said to be
strictly outer when the relative commutant of the factor in the
crossed product is trivial. We show that
all locally compact quantum groups can act strictly outerly on a free
Araki-Woods factor and that all locally compact
groups can act strictly outerly on the hyperfinite II$_1$ factor. We
define a kind of Connes' $T$ invariant for locally compact quantum groups and
prove a link with the possibility of acting strictly outerly on a
factor with a given $T$ invariant. Necessary and sufficient conditions
for the existence of strictly outer actions of compact Kac algebras on
the hyperfinite II$_1$ factor are given.
\end{abstract}

\maketitle

\section{Introduction}

The most natural appearance of a group is as a symmetry group of a
space. Several types of groups are simply defined by their action on a
space. More generally, actions of locally compact (l.c.)\ groups on
quantum spaces, i.e.\ von Neumann algebras, attracted a lot of
attention. Therefore, it is a very natural idea to consider also
actions of l.c.\ quantum groups on von Neumann algebras.

The history of the theory of l.c.\ quantum groups dates back to the
1960's, with the important work of Kac, Kac \& Vainerman and Enock \&
Schwartz, see \cite{E-S} for an overview. Building on important
contributions of Baaj \& Skandalis \cite{BS}, Woronowicz
\cite{Wor1,Wor2} and Van Daele \cite{VD}, the general theory of l.c.\
quantum groups was developed by Kustermans and the author
\cite{KV1,KV2}. The general theory of actions of l.c.\ quantum groups
on von Neumann algebras was developed in \cite{V}.

It has been shown by Enock and Nest \cite{EN,E} that quantum group
symmetries appear in the study of irreducible inclusions of factors,
of depth~2 and
infinite index. If $N_0 \subset N_1$ is such an inclusion and satisfies a technical regularity
condition, they consider the Jones tower $N_0 \subset N_1 \subset N_2
\subset \cdots$. On the relative commutant $N_2 \cap N_0'$, Enock and
Nest construct a l.c.\ quantum group structure. They construct a
strictly outer action of this l.c.\ quantum group on $N_1$ such that
$N_0$ is the fixed point algebra. This link between quantum groups and
irreducible inclusions of depth~2 goes back to Ocneanu and was also
generalized to reducible inclusions, yielding quantum groupoids, see
e.g.\ \cite{NV}.

It is now an obvious question to study which l.c.\ quantum groups can
appear in Enock and Nest's construction? Equivalently (see \cite{V}),
the question is which l.c.\ quantum groups can act \emph{strictly outerly} on
a factor, i.e. \emph{such that the relative commutant of the factor in the
crossed product is trivial.}

Even for l.c.\ groups, this question is not so easy. Blattner \cite{Bla} constructed for any l.c.\ group $G$, an action $(\al_g)$ of $G$ on the
hyperfinite II$_1$ factor $\cR$ and showed that none of the $\al_g$ is an inner automorphism (for $g \neq e$). This is a necessary condition in order
to get a trivial relative commutant of the factor in the crossed product, but far from being a sufficient condition.

For compact Lie groups, the question is more easy to answer: if $g \mapsto u_g$ is a finite-dimensional representation, say on $\C^n$, such that $g
\mapsto \Ad u_g$ is faithful, Wassermann (\cite{Was}, page 212) considered the diagonal action $\Ad u_g$ on the infinite tensor product $\cR =
\bigotimes_{n=1}^\infty \M_n(\C)$ and showed that this action is minimal. This means that the relative commutant of the fixed point algebra in $\cR$
is trivial and implies, in particular, that the action is strictly outer.

For certain actions of $\R$ on the hyperfinite II$_1$ factor,
Kawahigashi proved in
Proposition 3.2 of \cite{Kaw} the strict outerness. It might be
possible to use Kawahigashi's techniques to prove that Blattner's
action is strictly outer for every infinite abelian l.c.\
group.

Our contribution to the case of l.c.\ groups consists in proving the
following result: if $(\al_g)$ is an action of a l.c.\ group $G$ on a
factor $N$ leaving invariant a faithful normal state and such that $g \mapsto \al_g$ is
faithful, than the diagonal action of $G$ on the infinite tensor
product of copies of $N$ (w.r.t.\ $\om$) is strictly outer. As a
corollary, we obtain that all l.c.\ groups can act strictly outerly on
the hyperfinite II$_1$ factor. In fact, a variant of Blattner's
action, taking infinitely many copies of the regular representation,
is of this form and hence, strictly outer. Further, we give a more
geometrical construction of a strictly outer action of any linear group
on the hyperfinite II$_1$ factor.

We turn now to the quantum setting. For finite quantum groups,
Yamanouchi \cite{Y1} constructed a minimal (and hence, strictly outer)
action on the hyperfinite II$_1$ factor. Even to obtain strictly outer
actions of compact quantum groups is not so easy. For instance, it was shown recently by Izumi
\cite{Iz} that an infinite tensor product action in the style of
Wassermann, fails to be strictly outer for the quantum $SU_q(2)$
group.
Strictly outer actions of arbitrary compact quantum groups were constructed by Ueda \cite{U1}, by
taking the free product of the action of the quantum group on itself
by translation, and a trivial action. Of course, using the free
product, the obtained actions are no longer on hyperfinite factors.

Below, we prove that there exists a III$_1$ factor -- it is a free
Araki-Woods factor in the sense of Shlyakhtenko
\cite{shl} -- on which every l.c.\ quantum group can act strictly outerly.
In fact, the key point in the argument of Ueda is the
following: because the Haar measure of a compact quantum group is a
state and invariant under the comultiplication, it is possible to make
a free product w.r.t.\ this invariant state. In the non-compact
case, we can no longer use the Haar measure, because it is infinite
and the free product construction w.r.t.\ weights does not work. So,
we have to produce first, in a different way, an action of an arbitrary l.c.\
quantum group on a von Neumann algebra leaving invariant a faithful
state. Afterwards, we can apply the idea of Ueda and take the free
product with a trivial action.

Having strictly outer actions of an arbitrary l.c.\ quantum group on a
factor, it is a natural idea to ask when such a strictly outer action
can be found on certain types of factors. The general construction
gives a strictly outer action on a type III$_1$ factor.
Of course, a strictly outer
action cannot exist on a factor of type I. Below, we give necessary
and sufficient conditions for a l.c.\ quantum group being able to act
strictly outerly on a factor of type II$_1$, II$_\infty$, III$_\sla$
for $0 < \sla < 1$ or certain III$_0$ factors. In fact, we introduce a kind of Connes $T$
invariant $T(M,\de)$ for a l.c.\ quantum group $(M,\de)$ (which is not
the $T$ invariant of the von Neumann algebra $M$). We prove that if a
l.c.\ quantum group $(M,\de)$ acts strictly outerly on a factor $N$,
then $T(N) \subset T(M,\de)$.

As we already remarked, the free Araki-Woods factors on which we find strictly outer actions of arbitrary l.c.\ quantum groups, are not injective.
Since every l.c.\ group can act strictly outerly on an injective factor, it is a natural question to study which l.c.\ \emph{quantum} groups can act
strictly outerly on an injective factor. We prove that a necessary condition is the co-amenability of the l.c.\ quantum group. The converse is
however far from clear (and probably, false), but we prove that a compact Kac algebra can act strictly outerly on an injective factor if and only if
it is co-amenable. We further show that the bicrossed product
l.c.\ quantum groups \cite{VV,BSV} (which need not be Kac algebras)
can act strictly outerly on an injective factor if and only if they
are co-amenable. Using results of Ueda \cite{U2}, we show that if the compact quantum group
$SU_q(2)$ acts strictly outerly on an injective factor, then there exists an irreducible subfactor of the hyperfinite II$_1$ factor with index
$(q+q^{-1})^2$. Since most important unpublished work of Popa makes this last
statement highly improbable, we are convinced that not all the
$SU_q(2)$ can act strictly outerly on an injective factor.

Observe that Banica (\cite{Ban}, Section 4) has shown that a
discrete Kac algebra with a faithful finite dimensional
corepresentation (see Definition \ref{def.faithfulcorep}) can act strictly outerly on the hyperfinite II$_1$
factor. We generalize this result and prove that a discrete Kac
algebra with an faithful corepresentation in the
hyperfinite II$_1$ factor $\cR$ can act strictly outerly on
$\cR$. Such discrete Kac algebras need not be amenable.

{\bf Acknowledgment } The author is highly grateful to Georges Skandalis
for all the discussions and his essential contribution to the results
of this paper.

\section{Preliminaries}
For simplicity, we will assume that all Hilbert spaces are separable
and all von Neumann algebras have separable predual. So, we only work
with second countable locally compact quantum groups. These
assumptions are not essential but give us lighter statements of our results.

We study actions of locally compact (l.c.) groups and quantum groups
on von Neumann algebras. An action of $G$ on a von Neumann algebra is
a morphism $G \recht \operatorname{Aut} G : g \mapsto \al_g$ such that
$g \mapsto \al_g(x)$ is strongly$^*$ continuous for every $x \in
N$. We then define the crossed product
$$G \kruisje{\al} N := \bigl( \al(N) \cup \cL(G) \ot 1
\bigr)^{\prime\prime} \subset \B(L^2(G)) \ot N \; ,$$
where $\al : N \recht L^\infty(G) \ot N : (\al(x))(g) =
\al_{g^{-1}}(x)$ and $\cL(G)$ is the group von Neumann algebra
generated by the \emph{left regular representation} on $L^2(G)$ for
the left Haar measure, defined by $(\sla_g \xi)(h) = \xi(g^{-1} h)$.

\begin{definition}
An action $(\al_g)$ of a l.c.\ group $G$ on a factor $N$ is called
\emph{strictly outer} if the relative commutant of $N$ in the crossed
product is trivial, i.e.
$$G \kruisje{\al} N \cap \al(N)' = \C \; .$$
 \end{definition}
It is obvious that a strictly outer action is outer, i.e.\ if $g \neq
e$, then $\al_g$ is not an inner automorphism. But, if $G$ is not a
discrete group, not every outer action is strictly outer. For
instance, if $N$ is a III$_0$ factor with trivial $T$ invariant, the
modular automorphism group $(\si^\te_t)$ of an n.s.f.\ weight $\te$ on
$N$ is outer ($\si^\te_t$ is not inner for $t \neq 0$), but not
strictly outer: the crossed product is even not a factor, but its
center is the space of the flow of weights.

We generalize these concepts to the world of l.c.\ quantum groups. The
general theory of l.c.\ quantum groups was developed by Kustermans and
the author in \cite{KV1,KV2}. We recall some of the basic definitions
and properties.

\begin{definition}
A pair $(M,\de)$ is called a (von Neumann algebraic) l.c.\ quantum group when
\begin{itemize}
\item $M$ is a von Neumann algebra and $\de : M \recht M \ot M$ is
a normal and unital $*$-homomorphism satisfying the coassociativity relation : $(\de \ot \io)\de = (\io \ot \de)\de$.
\item There exist normal semi-finite faithful (n.s.f.) weights $\varphi$ and $\psi$ on $M$ such that
\begin{itemize}
\item $\varphi$ is left invariant in the sense that $\varphi \bigl( (\om \ot
\io)\de(x) \bigr) = \varphi(x) \om(1)$ for all $x \in \Mfi^+$ and $\om \in M_*^+$,
\item $\psi$ is right invariant in the sense that $\psi \bigl( (\io \ot
\om)\de(x) \bigr) = \psi(x) \om(1)$ for all $x \in \Mpsi^+$ and $\om \in M_*^+$.
\end{itemize}
\end{itemize}
\end{definition}
Ordinary l.c.\ groups appear in this theory in the form of $M =
L^\infty(G)$ and $(\de(f))(p,q) = f(pq)$ for $f \in L^\infty(G)$. All
l.c.\ quantum groups whose von Neumann algebra $M$ is commutative are
of this form.

Fix a l.c.\ quantum group $(M,\de)$.

We first define the analogue of the left regular representation. As
usual, we write $\Nfi = \{ x \in M \mid \vfi(x^*x) < \infty \}$.
Represent $M$ in the GNS-construction of $\vfi$ with GNS-map $\la : \Nfi \recht H$. We define a unitary $W$ on $H \ot H$ by
$$W^* (\Lambda(a) \ot \Lambda(b)) = (\Lambda \ot \Lambda)(\de(b)(a \ot 1)) \quad\text{for all}\; a,b \in N_{\phi}
\; .$$ Here, $\Lambda \ot \Lambda$ denotes the canonical GNS-map for the tensor product weight $\varphi \ot \varphi$. One proves that $W$ satisfies
the pentagonal equation: $W_{12} W_{13} W_{23} = W_{23} W_{12}$, and we say that $W$ is a \emph{multiplicative unitary}. It is the \emph{left regular
corepresentation}. The von Neumann algebra $M$ is the strong closure of the algebra $\{ (\io \ot \om)(W) \mid \om \in \B(H)_* \}$ and $\de(x) = W^* (1
\ot x) W$, for all $x \in M$. The {\it von Neumann algebraic} quantum
group $(M,\de)$ has an underlying {\it C$^*$-algebraic} quantum group
$(A,\de)$, where $A$ is the norm closure of $\{ (\io \ot \om)(W) \mid
\om \in \B(H)_* \}$. If $M=L^\infty(G)$, we have $A = C_0(G)$, the
continuous functions on $G$ vanishing at infinity.

Next, the l.c.\ quantum group $(M,\de)$ has an {\it antipode} $S$, which is the unique $\si$-strong$^*$ closed linear map
from $M$ to $M$ satisfying $(\io \ot \om)(W) \in D(S)$ for all $\om \in \B(H)_*$ and $$S\bigl((\io \ot \om)(W)\bigr) = (\io \ot \om)(W^*)$$ and such that the
elements $(\io \ot \om)(W)$ form a $\si$-strong$^*$ core for $S$. The antipode $S$ has a polar decomposition $S = R \tau_{-i/2}$ where $R$ is an
anti-automorphism of $M$ and $(\tau_t)$ is a strongly continuous one-parameter group of automorphisms of $M$. We call $R$ the \emph{unitary antipode}
and $(\tau_t)$ the \emph{scaling group} of $(M,\de)$.

The \emph{dual l.c.\ quantum group} $(\Mh,\deh)$ is defined in \cite{KV1}, Section~8. Its von Neumann algebra $\Mh$ is the strong closure of the
algebra $\{(\om \ot \io)(W) \mid \om \in \B(H)_* \}$ and the comultiplication is given by $\deh(x) = \Sigma W (x \ot 1) W^* \Sigma$ for all $x \in
\Mh$,
where $\Si$ is the flip map on $H \ot H$. On $\Mh$ there exists a canonical left invariant weight $\vfih$
and the associated multiplicative unitary is given by $\Si W^*
\Si$. We have again an underlying C$^*$-algebraic quantum group
$(\Ah,\deh)$ where $\Ah$ is the norm closure of $\{(\om \ot \io)(W) \mid \om \in \B(H)_* \}$.

Since $(\Mh,\deh)$ is again a l.c.\ quantum group, we can introduce
the antipode $\hat{S}$, the unitary antipode $\hat{R}$ and the scaling group
$(\hat{\tau}_t)$ exactly as we did it for $(M,\de)$. Observe that
$$\Sh \bigl((\om \ot \io)(W^*)\bigr) = (\om \ot \io)(W) \; .$$

The modular conjugations of the weights $\varphi$ and $\hat\varphi$
will be denoted by $J$ and $\Jh$ respectively. Then it is worthwhile to mention that $$R(x) = \Jh x^* \Jh \quad\text{for all} \; x \in M
\qquad\text{and}\qquad \Rh(y) = J y^* J \quad\text{for all}\; y \in
\Mh \; .$$ The modular operators of the weights $\vfi$ and $\vfih$ are
denoted by $\nab$ and $\nabh$. We mention that
$$\tau_t(x) = \nabh^{it} x \nabh^{-it} \quad\text{for all} \; x \in M
\qquad\text{and}\qquad \tauh_t(x) = \nab^{it} y \nab^{-it} \quad\text{for all}\; y \in
\Mh \; .$$

\begin{definition} \label{def.faithfulcorep}
A unitary $U \in M \ot \B(K)$ is called a {\it corepresentation} of a l.c.\ quantum group $(M,\de)$ on the
Hilbert space $K$, if $(\de \ot \io)(U) = U_{13}
U_{23}$.

We say that $U$ is {\it faithful} if
$$\{(\io \ot \mu)(U) \mid \mu \in \B(K)_* \}\dpr = M \; .$$
\end{definition}

We observe that a corepresentation $U$ is, in a sense,
automatically continuous: we have $U \in \M(A \ot \cK(K))$, where
$\cK(K)$ denotes the C$^*$-algebra of compact operators.

Taking into account all unitary corepresentations, we can
define a {\it universal C$^*$-algebraic dual} $(\Ahu,\dehu)$. The non-degenerate $^*$-representations of the
C$^*$-algebra $\Ahu$, $\rho : \Ahu \recht \B(K)$ are in one-to-one
correspondence with the corepresentations of $U$ of $(M,\de)$ through
the formula $U = (\io \ot \rho)(\cW)$, where $\cW \in \M(A \ot \Ahu)$
is the so-called {\it universal corepresentation}. See \cite{JK} for
details. If $M= L^\infty(G)$, the universal C$^*$-algebraic dual is
$C^*(G)$, the full group C$^*$-algebra.

Of course, since $(M,\de)$ has a universal dual $(\Ahu,\dehu)$, also
$(\Mh,\deh)$ has a universal dual $(\Au,\deu)$ and we have a universal
corepresentation $\cWh \in \M(\Ah \ot \Au)$.

We use \cite{V} as a reference for actions of l.c.\ quantum groups, but we recall the necessary elements of the theory. Let $(M,\de)$ be a l.c.\
quantum group and $N$ a von Neumann algebra. A faithful, normal $^*$-homomorphism $\al : N \recht M \ot N$ is called a {\it (left) action} of
$(M,\de)$ on $N$ if $$(\de \ot \io)\al = (\io \ot \al)\al \; .$$ One can define the crossed product as
$$M \kruisje{\al} N = (\al(N) \cup \Mh \ot 1)\dpr \subset \B(H) \ot N \; .$$

If $\al : N \recht M \ot N$ is an action, the \emph{fixed point algebra} is denoted by $N^\al$ and defined as the von Neumann subalgebra of elements
$x \in N$ satisfying $\al(x) = 1 \ot x$.

We repeat Definition 6.1 of \cite{V}.

\begin{definition} \label{def.faithmin}
An action $\al$ of $(M,\de)$ on $N$ is called \emph{faithful} when
$$\{ (\io \ot \om)\al(x) \mid \om \in N_*, x \in N \}^{\prime\prime} =
M \; .$$
The action is called \emph{minimal} when it is faithful and
$N \cap (N^\al)' = \C$.
\end{definition}

We finally define strictly outer actions of l.c.\ quantum groups.

\begin{definition}
An action of a l.c.\ quantum group
$(M,\de)$ on a factor $N$ is called \emph{strictly outer} if
the relative commutant of this factor in the crossed
product is trivial, i.e.
$$M \kruisje{\al} N \cap \al(N)' = \C \; .$$
\end{definition}

From Proposition 6.2 in \cite{V}, we know that every minimal action is strictly outer. An integrable strictly outer action is minimal. In particular,
\emph{for compact quantum groups} the notions of minimal and strictly outer actions coincide.

Observe also that an action of an ordinary l.c.\ group
$G$ on a von Neumann algebra $N$ is faithful in the sense of
Definition \ref{def.faithmin} if and only if the morphism $g \mapsto
\al_g$ is faithful.

An important technical tool in the study of actions of l.c.\ quantum groups on von Neumann algebras is the Radon-Nikodym derivative of an n.s.f.\
weight under the action, as introduced by Yamanouchi \cite{Y} (see the appendix of \cite{BV} for an easy approach). If $\al : N \recht M \ot N$ is an
action of $(M,\de)$ on $N$ and if $\te$ is an n.s.f.\ weight on $N$, we can consider unitaries $D_t:=[D \te \circ \al : D \te ]_t \in M \ot N$,
satisfying the following properties:
$$(\de \ot \io)(D_t) = (\io \ot \al)(D_t) \; (1 \ot D_t) \qquad\text{and}\qquad D_{t+s} = D_t \; (\tau_t \ot \si^\te_t)(D_s) \; .$$
If we represent the crossed product $M \kruisje{\al} N$ as above, as a subalgebra of $\B(H) \ot N$, we can define a dual n.s.f.\ weight $\tilde{\te}$
on $M \kruisje{\al} N$ with a modular operator $\nabtil$ such that $D_t = \nabtil^{it} (\nabh^{-it} \ot \nab_\te^{-it})$, where $\nabh$ is the
modular operator of the left invariant weight on the dual $(\Mh,\deh)$ and $\nab_\te$ is the modular operator of $\te$.

\begin{definition}
Let $\al : N \recht M \ot N$ be an action of $(M,\de)$ on $N$ and let
$\te$ be an n.s.f.\ weight on $N$.
We say that $\te$ is $\rho$-invariant, when $\rho$ is a strictly
positive, self-adjoint operator affiliated with $M$, satisfying
$$\de(\rho)=\rho \ot \rho \quad\text{and}\quad \te \bigl(
(\om_{\xi,\xi} \ot \io)\al(x) \bigr) = \| \rho^{1/2} \xi \|^2 \;
\te(x) \quad\text{for}\quad x \in \cM_\te^+,\xi \in D(\rho^{1/2}) \;
.$$
We say that a state $\om$ on $N$ is invariant, when $(\io \ot
\om)\al(x) = \om(x) \; 1$ for all $x \in N$.
\end{definition}

We finish this section of preliminaries with the following alternative
characterization of strict outerness.

\begin{proposition} \label{prop.alternativeouter}
Let $\al : N \recht M \ot N$ be an action of a l.c.\ quantum group on
a factor $N$. Then, $\al$ is strictly outer if and only if
$$\B(H) \ot N \cap \al(N)' = M' \ot 1 \; .$$
In particular, a strictly outer action is faithful.
\end{proposition}
\begin{proof}
Suppose first that $\B(H) \ot N \cap \al(N)' = M' \ot 1$. In order to
prove that $\al$ is strictly outer, it suffices to show that $M
\kruisje{\al} N \cap M' \ot 1 = \C 1$. If $V \in \Mh' \ot M$ denotes
the right regular corepresentation of $(M,\de)$ (see \cite{KV2}), it is
easy to check that $(\io \ot \al)(z) = V_{12} z_{13} V_{12}^*$ for all
$z \in M \kruisje{\al} N$. So, if $a \in M'$ and $a \ot 1 \in M
\kruisje{\al} N$, we get $a \in M' \cap \Mh = \C 1$ and we are done.

Suppose conversely that $\al$ is strictly outer and let $z \in \B(H)
\ot N \cap \al(N)'$. Then,
$$W_{12} (\io \ot \al)(z) W_{12}^* \in \B(H) \ot (M \kruisje{\al} N)
\cap (1 \ot \al(N))' = \B(H) \ot 1 \ot 1 \; .$$
So, we can take $a \in \B(H)$ such that $(\io \ot \al)(z) = W^*(a \ot
1)W \ot 1$. The left hand side belongs to $\B(H) \ot M \ot N$, while
the right hand side belongs to $\B(H) \ot \Mh \ot 1$. Since $M \cap
\Mh = \C$, the right
hand side belongs to $\B(H) \ot 1 \ot 1$ and we can take $b \in \B(H)$
such that $z = b \ot 1 = W^*(a \ot 1) W$. If we apply $\io \ot \deh$,
we conclude that
$$b \ot 1 \ot 1 = (\io \ot \deh)(W^*(a \ot 1) W) = W^*_{12} W^*_{13}
(a \ot 1 \ot 1) W_{13} W_{12} = W^*_{12} (b \ot 1 \ot 1) W_{12} \; .$$
It follows that $b \in M'$ and hence, $z \in M' \ot 1$. This proves
the first part of the proposition.

If $a$ commutes with $(\io \ot \mu)\al(x)$ for all $\mu \in N_*$ and
$x \in N$, it is clear that $a \ot 1 \in \B(H) \ot N \cap
\al(N)'$. Hence, $a \in M'$. The faithfulness of $\al$ follows
immediately.
\end{proof}

\section{Strictly outer actions of locally compact quantum groups} \label{sec.outeractions}

We shall prove,
in Theorem \ref{thm.construction}.a) below, that there exists a type
III$_1$ free Araki-Woods factor (in the sense of Shlyakhtenko
\cite{shl}) on which every l.c.\ quantum group acts strictly
outerly. We introduce a T invariant for l.c.\ quantum groups
(Definition \ref{def.Tinvariant}) and determine when a given l.c.\
quantum group can act strictly outerly on a factor of a given type
(Theorems \ref{thm.restrictions} and \ref{thm.construction}).

Let us recall Shlyakhtenko's construction of these free Araki-Woods
factors. Let $T$ be an involution on a Hilbert space $K$. This means
that $T$ is a densely defined, closed, injective, anti-linear operator
on $K$ satisfying $T^{-1}=T$. If $T = \cJ Q^{1/2}$ is its polar
decomposition, $\cJ$ is an anti-unitary operator and $Q$ is a strictly
positive, self-adjoint operator satisfying $\cJ Q \cJ = Q^{-1}$. We
also have $\cJ Q^{it} \cJ = Q^{it}$. Define
$$K_\cJ := \{ \xi \in H \mid \cJ \xi= \xi \} \; .$$
Then, $K_\cJ$ is a real Hilbert space and the restriction of $Q^{it}$
gives a one-parameter group of orthogonal transformations of $K_\cJ$.

Consider the full Fock space
$$\cF(K) := \C \Om \oplus \bigoplus_{n=1}^\infty K^{\ot (n)} \; .$$
Shlyakhtenko's free Araki-Woods von Neumann algebra \cite{shl} can be defined as
$$\free{K_\cJ,Q^{it}} := \{ s(\xi) \mid \xi \in D(T) \}\dpr \subset \B(\cF(K))
\quad\text{with}\quad s(\xi) = \ell(\xi)+\ell(T\xi)^* \; ,$$
where $\ell(\xi)$ denotes the left creation operator in $\B(\cF(K))$
whenever $\xi \in K$. The vector state $\om_\Omega$ is faithful on
$\free{K_\cJ,Q^{it}}$ and is called the free quasi-free state.

Instead of starting with an involution $T$, we can of course start
with a real Hilbert space with a one-parameter group of orthogonal
operators and complexify to obtain $T$.

We show now that any l.c.\ quantum group can act in a natural way on a free Araki-Woods
von Neumann algebra, whenever we have a corepresentation of $(M,\de)$
on a Hilbert space $K$ and a compatible involution $T$ on $K$.

\begin{proposition} \label{prop.action}
Let $U \in M \ot \B(K)$ be a corepresentation of a l.c.\ quantum group
$(M,\de)$ on a Hilbert space $K$. Suppose that there exists an
involution $T$ on $K$ satisfying
\begin{equation} \label{eq.crucial}
(\mu \ot \io)(U^*) T  \subset T (\overline{\mu} \ot \io)(U^*)
\end{equation}
for all $\mu \in M_*$. Let $T = \cJ Q^{1/2}$ be the polar
decomposition. Denote by $\cF(U)$ the amplified corepresentation on
the full Fock space $\cF(K)$ defined by
$$\cF(U) = \bigoplus_{n=0}^\infty U^{(n)} \; ,$$
where $U^{(0)} = 1 \in M \ot \B(\C \Omega)$ and $U^{(n)} = U_{1,n+1}
\ldots U_{12} \in M \ot \B(K^{\ot (n)})$ for $n \geq 1$. Let $N:=
\free{K_\cJ,Q^{it}}$ be the free Araki-Woods factor corresponding to $T$. Defining
$$\al : N \recht M \ot N : \al(z) = \cF(U)^* (1 \ot z) \cF(U) \; ,$$
we get a well-defined action of $(M,\de)$ on $N$ leaving invariant the
free quasi-free state.

If the corepresentation $U$ is faithful, the action
$\al$ is faithful.
\end{proposition}
\begin{proof}
It is obvious that $\cF(U)$ is a corepresentation of $(M,\de)$ on the
full Fock space $\cF(K)$.
We observe that
$$(\mu \ot \io)\bigl( \cF(U)^* (1 \ot \el(\xi)) \cF(U) \bigr) =
\el\bigl( (\mu \ot \io)(U^*) \xi \bigr)\; .$$ Using
Equation~\eqref{eq.crucial}, we get that
$$(\mu \ot \io)\bigl( \cF(U)^* (1 \ot s(\xi)) \cF(U) \bigr) =
s\bigl( (\mu \ot \io)(U^*) \xi \bigr) \; ,$$
for all $\mu \in M_*$ and $\xi \in D(T)$. Hence, the definition of
$\al$ in the statement of the proposition yields an action of
$(M,\de)$ on the von Neumann algebra $\free{K_\cJ,Q^{it}}$. It is obvious that the
free quasi-free state is invariant under the action
$\al$.

Further, we observe that $(\io \ot \om_{\Om,\eta})(\al(s(\xi))) = (\io \ot
\om_{\xi,\eta})(U^*)$ for $\xi \in D(T)$ and $\eta \in H
\hookrightarrow \cF(H)$. So, the action $\al$ is faithful if the
corepresentation $U$ is faithful.
\end{proof}

\begin{remark}
The assumption on the existence of an involution $T$ satisfying
Equation \eqref{eq.crucial} is crucial for the construction of the
action $\al$. It has a very natural interpretation, as follows. Let
$\rho : \Ahu \recht \B(K)$ be the representation of the universal
C$^*$-algebraic dual $(\Ahu,\dehu)$ corresponding to $U$ through the
formula $(\io \ot \rho)(\cW) = U$ (see Preliminaries). From \cite{JK},
we know that $(\Ahu,\dehu)$ has an antipode $\Shu$ satisfying
$\Shu\bigl( (\mu \ot \io)(\cW^*) \bigr) = (\mu \ot \io)(\cW)$ for all
$\mu \in A^*$. This antipode has a polar decomposition $\Shu = \Rhu
\tauhu_{-i/2}$. It is not hard to check that Equation
\eqref{eq.crucial} is then equivalent with the equations
$$\rho(\Rhu(x)) = \cJ \rho(x)^* \cJ \quad\text{and}\quad
\rho(\tauhu_t(x)) = Q^{it} \rho(x) Q^{-it} \;\; \text{for all}\;\; x
\in \Ahu \; \text{and}\; t \in \R \; .$$
If we define the real C$^*$-algebra
$$\Ah\uni\real := \{ x \in \Ah\uni \mid \Rhu(x) = x^* \}$$
and the real Hilbert space $K_\cJ$ as before, we see that
$\rho(\Ah\uni\real)$ leaves $K_\cJ$ invariant. Further, the
one-parameter group $(\tauhu_t)$ restricts to a one-parameter group of
$\Ah\uni\real$. So, a compatible pair of a corepresentation $U$ of
$(M,\de)$ on a Hilbert space $K$ and an involution $T=\cJ Q^{1/2}$ is
the same thing as a representation $\rho$ of the real C$^*$-algebra
$\Ah\uni\real$ on the real Hilbert space $K_\cJ$ satisfying
$\rho(\tauhu_t(x)) = Q^{it} \rho(x) Q^{-it}$ for all $x \in
\Ah\uni\real$ and $t \in \R$.

Finally, the involution $T$ exists in an important case.
If we
take the left regular corepresentation $W \in \B(H \ot H)$, we know
that Equation \eqref{eq.crucial} holds for $W$ with $T =J
\nab^{1/2}$. This is the case because $J$ implements the unitary
antipode $\Rh$ of $(\Mh,\deh)$, while $\nab^{it}$ implements the
scaling group of $(\Mh,\deh)$ and $\Sh((\om \ot \io)(W^*)) = (\om \ot \io)(W)$. We
refer to the section of preliminaries for the definition of $W$, $J$
and $\nab$.
\end{remark}

\begin{corollary} \label{cor.faithaction}
Every l.c.\ quantum group can act faithfully on
$\free{L^2(\R,K_\R),\sla_t}$ leaving invariant the free quasi-free
state. Here, $K_\R$ is the separable, infinite-dimensional real
Hilbert space and $\sla_t$ acts by left translations on $\R$.

Every l.c.\ quantum group with trivial scaling group $\tau_t = \io$
can act faithfully on the free group factor $\cL(\F_\infty)$ leaving
invariant the trace.
\end{corollary}
\begin{proof}
Consider the Hilbert space $L^2(\R,H)$ and define $(\cJ \xi)(x) = J
\xi(x)$, $(Q^{it} \xi)(x) = \nab^{it} \xi(x-t)$ for $\xi \in
L^2(\R,H)$. Then, $T:= \cJ Q^{1/2}$ is an involution. Identifying
$L^2(\R,H)$ with $L^2(\R) \ot H$, we define a corepresentation
$U:=W_{13}$ of $(M,\de)$ on $L^2(\R,H)$. We already remarked above
that Equation \eqref{eq.crucial} is satisfied for $W$ and
$J\nab^{1/2}$. Hence, Equation \eqref{eq.crucial} is satisfied for our
$U$ and $T$. From Proposition \ref{prop.action}, we get a faithful
action of $(M,\de)$ on the free Araki-Woods factor
$\free{L^2(\R,H)_\cJ,Q^{it}}$ leaving invariant the free quasi-free
state.

Using the unitary $(Z\xi)(x) = \nab^{-ix} \xi(x)$ on $L^2(\R,H)$, it
is clear that
$$(L^2(\R,H)_\cJ,Q^{it}) \cong (L^2(\R,H_\cJ),\sla_t) \; .$$
So, the first part of the corollary is proven. If $M$ happens to be finite dimensional,
we take above a direct sum of infinitely many copies of $W$ and assure
as such that we can always assume $K_\R$ to be infinite dimensional.

To prove the second part, suppose that $\tau_t = \io$ for all
$t$. Then, $\tauh_t=\io$ as well and $(\mu \ot \io)(W^*) J = J
(\overline{\mu} \ot \io)(W^*)$ for $\mu \in M_*$. By proposition
\ref{prop.action}, we get an action of $(M,\de)$ on $\free{H_J,\io}$
leaving invariant the trace. As above, we can assume that $H$ is
infinite dimensional and then, $\free{H_J,\io} \cong \cL(\F_\infty)$.
\end{proof}

The following observation is due to Y.\ Ueda \cite{U1}.
If $\al$ is a faithful action of $(M,\de)$ on a von Neumann algebra
$N_0$ leaving invariant a faithful state $\om_0$ and if $N_1$
is an arbitrary von Neumann algebra with faithful state $\om_1$, it is
easy to extend $\al$ to an action $\be$ on the free product (see
\cite{voi}) $(N,\om):= (N_0,\om_0)
\star (N_1,\om_1)$ acting trivially on $N_1$. If now, $N \cap N_1'
= \C$, the action $\be$ is minimal and hence, strictly outer. It is
easy to give examples such that $N \cap N_1' = \C$. Indeed, we can
take $(N_1,\om_1) = (A,\eta) \star (B,\mu)$ such that the centralizer
$A^\eta$ contains a non-trivial group of orthogonal unitaries and such
that the centralizer $B^\mu$ contains a group of orthogonal unitaries
with at least three elements. It follows from a result of Barnett (Theorem 11 in \cite{Bar})
that
$$N \cap N_1' = (N \star A) \star B \cap (A \star B)' = \C \;
.$$
Moreover, $N$ is a full factor.

We conclude that every l.c.\ quantum group can act strictly outerly on
a full factor.

In Theorem \ref{thm.construction}.a) below, we want to construct such a strictly outer action on the
canonical free Araki-Woods factor $\free{L^2(\R,K_\R),\sla_t}$ and for
this, we will need a generalization of Barnett's result, see
Lemma \ref{lemma.barnett}.

As we know now that every l.c.\ quantum group can act strictly outerly on a factor, it is a natural question to study the possibility of acting
strictly outerly on factors with certain given invariants. We introduce some kind of Connes' $T$ invariant of a l.c.\ quantum group and show how it
is related to the possibility of acting strictly outerly on a factor with a certain $T$ invariant.

\begin{definition} \label{def.Tinvariant}
Let $(M,\de)$ be a l.c.\ quantum group and $(\tau_t)$ its scaling group. We define
$$T(M,\de) := \{ t \in \R \mid \; \text{There exists a unitary} \;\; u \in M \;\; \text{satisfying} \; \; \tau_t = \Ad u \;\;\text{and}\;\; \de(u) =
u \ot u \} \; .$$
\end{definition}
Remark that, exactly as in the case of the usual $T$ invariant of von
Neumann algebras defined by Connes in \cite{Con}, we have that $T(M,\de)$ is
a subgroup of $\R$.

Since $\de \tau_t = (\tau_t \ot \tau_t)\de$, it is natural to
consider, in the definition of $T(M,\de)$ that $\tau_t$ is implemented
by a unitary satisfying $\de(u) = u \ot u$.

We denote by $T(N)$ the $T$ invariant of a von Neumann algebra $N$.

In the following result, we prove that a small $T$ invariant
$T(M,\de)$ makes it impossible to act on certain factors. More
precisely, we have the following.

\begin{theorem} \label{thm.restrictions}
Let $(M,\de)$ be a l.c.\ quantum group.
\begin{enumerate}\renewcommand{\theenumi}{\alph{enumi}}
\item If $(M,\de)$ acts strictly outerly on a factor $N$, then $T(N) \subset
T(M,\de)$. \label{en.a}
\item If $(M,\de)$ acts strictly outerly on a II$_\infty$ factor, there
exists a one-parameter group $\rho^{it} \in M$ of unitaries such that $\de(\rho^{it}) =
  \rho^{it} \ot \rho^{it}$ and $\tau_t = \Ad \rho^{it}$ for all
  $t \in \R$. \label{en.b}
\item If $(M,\de)$ acts strictly outerly on a II$_1$ factor, the scaling
group $(\tau_t)$ is trivial. \label{en.c}
\end{enumerate}
\end{theorem}
\begin{proof}
\ref{en.a}) Let $\al$ be a strictly outer action of $(M,\de)$ on a factor $N$ and $t_0 \in T(N)$. Take an n.s.f.\ weight $\te$ on $N$ such that
$\site_{t_0} = \io$. Then, we have a dual n.s.f.\ weight $\tetil$ on the crossed product $M \kruisje{\al} N$, with modular operator $\nabtil$. Denote
by $D_t := [D \te \circ \al : D \te]_t$ the Radon-Nikodym derivative and by $\nabte$ the modular operator of $\te$. As above, $\nabh$ denotes the
modular operator of the left invariant weight on the dual $(\Mh,\deh)$. We know that $\nabtil^{it} = D_t (\nabh^{it} \ot \nabte^{it})$. Take $x \in
N$. Because $D_t \in M \ot N$, we have $D_t (\nabh^{it} \ot 1) \in \B(H) \ot N$. Further, because $\nabte^{it_0}=1$ and $\sitil^\te_t \al = \al
\site_t$, we have
$$D_{t_0} (\nabh^{it_0} \ot 1) \; \al(x) = \nabtil^{it_0} \; \al(x) = \al(x) \; \nabtil^{it_0} = \al(x) \; D_{t_0} (\nabh^{it_0} \ot 1) \; .$$
From Proposition \ref{prop.alternativeouter} and
the strict outerness of $\al$, it follows that $\B(H) \ot N \cap \al(N)' = M' \ot 1$. So, because $D_{t_0} \in M \ot N$, there exists a unitary
$u \in M$ such that $D_{t_0} = u^* \ot 1$ and $u^* \nabh^{it_0} \in M'$. Because $\nabh^{it_0}$ implements $\tau_{t_0}$, we get $\tau_{t_0} = \Ad u$.
Finally, we know that $(\de \ot \io)(D_{t_0}) =$ \linebreak $(\io \ot \al)(D_{t_0}) (1 \ot D_{t_0})$ and so, $\de(u) = u \ot u$. We conclude that $t_0 \in
T(M,\de)$.

\ref{en.b}) In the reasoning above, we can take $\te=\Tr$, a trace on
$N$. We find that $D_t = u_t^* \ot 1$ for $u_t \in M$ and $t \in \R$, such that $\de(u_t) = u_t \ot u_t$ and $\tau_t
= \Ad u_t$. Because $D_{t+s} = D_t (\tau_t \ot \io)(D_s)$ and because any group-like unitary is invariant under the scaling group, we find that
$(u_t)$ is a one-parameter family of unitaries. So, we can write $u_t
= \rho^{it}$.

\ref{en.c}) Assume that, moreover, $\Tr$ is finite. Because $D_t = \rho^{-it} \ot 1$, it follows that $\Tr$ is
$\rho$-invariant. This implies that $\rho \; \Tr(1) = (\io \ot \Tr)\al(1) = 1 \; \Tr(1)$ and hence, $\rho=1$ and the scaling
group is trivial.
\end{proof}

From the constructive point of view, we prove the following. Recall
from \cite{shl}, Theorem 6.4 that for $0 < \sla < 1$, there exists a
unique free Araki-Woods factor of type III$_\sla$ given by
$\free{K_\cJ,Q^{it}}$ whenever the subgroup of $\R^*_+$ generated by
the spectrum of $Q$ is $\sla^\Z$.

\begin{theorem} \label{thm.construction}
Let $(M,\de)$ be a l.c.\ quantum group.
\begin{enumerate}\renewcommand{\theenumi}{\alph{enumi}}
\item Every l.c.\ quantum group can act strictly outerly on
  $\free{L^2(\R,K_\R),\sla_t}$ leaving invariant the free quasi-free
  state. \label{e.a}
\item If there exists a one-parameter group $\rho^{it} \in M$ of unitaries such that $\de(\rho^{it}) =
  \rho^{it} \ot \rho^{it}$ and $\tau_t = \Ad \rho^{it}$ for all
  $t \in \R$, $(M,\de)$ can act strictly outerly on $\cL(\F_\infty)
  \ot \B(K)$ leaving $\rho$-invariant the trace. \label{e.b}
\item If $t_0 \in T(M,\de)$, $(M,\de)$ can act strictly outerly on the
  unique free Araki-Woods factor $T_\sla$ of type III$_\sla$ with
  $0 < \sla < 1$ and $|t_0|=\frac{2\pi}{|\log \sla|}$. \label{e.c}
\item If $(\tau_t)$ is trivial, $(M,\de)$ can act strictly outerly on
  $\cL(\F_\infty)$ leaving invariant the trace. \label{e.d}
\item Let $\Gamma$ be a countable dense subgroup of $\R$. If $\Gamma \subset
  T(M,\de)$, $(M,\de)$ can act strictly outerly on a III$_0$ factor whose flow of weights is the
canonical action of $\R$ on the dual compact group $\Gammah$. \label{e.e}
\end{enumerate}
\end{theorem}
\begin{proof}
\ref{e.a}) From Corollary \ref{cor.faithaction}, we can take a
faithful action $\be$ of $(M,\de)$ on the free Araki-Woods factor $N_0 :=
\free{L^2(\R,K_\R),\sla_t}$ leaving invariant the free quasi-free
state $\om_0$. Define $N_1:=\free{L^2(\R,\R^3),\sla_t}$, with free
quasi-free state $\om_1$. Put $(N,\om) = (N_0,\om_0) \star
(N_1,\om_1)$. Then, $(N,\om) \cong (N_0,\om_0)$ by \cite{shl},
Theorem 2.11. Because $\om_0$ is invariant under $\be$, it is easy to
extend $\be$ to an action $\al$ of $(M,\de)$ on $N$, acting
trivially on $N_1$. It is clear that $\al$ is still
faithful. Moreover, if we write $A = \free{L^2(\R,\R),\sla_t}$, we
have $N = (N_0 \star A) \star A \star A$ and $A \star A \star A
\subset N^\al$ (all the free products are with respect to the free
quasi-free states). Combining Lemmas \ref{lemma.barnett} and
\ref{lemma.free} below, we get
that $N \cap (N^\al)' = \C$. So, $\al$ is a minimal action and
hence, a strictly outer action.

\ref{e.b}) To start with, take the above action $\al$ of $(M,\de)$ on
$N$ with invariant free quasi-free state $\om$. Define $C = N
\rtimes_{(\si^\om_t)} \R$ to be the continuous core of $N$. We
realize $C$ as the subalgebra of $N \ot \B(L^2(\R))$ generated by
$\si^\om(N)$ and $1 \ot \rho(\R)$, where $\si^\om(x) \in N \ot
L^\infty(\R)$ is defined by $\si^\om(x)(t) = \si_t^\om(x)$ and where
$\rho_t$ denotes the right translation by $t$.

Because the state $\om$ is invariant, we know from \cite{V}, proof
of Proposition 4.3, that $\be \si_t^{\om} = (\tau_t \ot \si_t^{\om})\be$. Defining the unitary $\cV
\in M \ot 1  \ot L^\infty(\R)$ by $\cV(s) = \rho^{is}  \ot 1$, we can define the action
\begin{align*}
& \ga : C \recht M \ot C : \ga(z) = \cV^* (\al \ot \io)(z) \cV
\quad\text{satisfying} \\ & \ga \si^{\om} = (\io \ot \si^{\om}) \al \; , \quad \ga(1 \ot \rho_t) = (\rho^{-it} \ot 1 \ot \rho_t) \; .
\end{align*}
It is clear that the dual weight of $\om$ is invariant under $\ga$ and
hence, the trace of $C$ will be $\rho$-invariant. Further,
$\si^\om(N^\al) \subset C^\ga$ and $\ga$ remains a faithful action. We
claim that $C \cap \si^\om(N^\al)' = \C$ and then, it follows that $\ga$
is a minimal and hence, strictly outer action. From Theorem 5.2 in
\cite{shl3}, we know that $C \cong \cL(\F_\infty) \ot \B(K)$. To
conclude the proof of item \ref{e.b}), it suffices hence to prove our
claim. Write again $N = (N_0 \star A) \star A \star A$. Take a
sequence $a_n$ of unitaries in the first copy of $A$ satisfying the
conclusions of Lemma \ref{lemma.free}. Take the same sequence in the second
and third copy of $A$ and denote them by $b_n$ and $c_n$. Let $x \in C
\cap \si^\om(N^\al)'$ and take $\mu \in \B(L^2(\R))_*$. Let $L$ be a
compact subset of $\R$ and write $p_L$ for the characteristic function
of $L$. Write $x_L = (1 \ot p_L)x(1 \ot p_L)$. Because $a_n \in
N^\al$, $x$ commutes with $\si^\om(a_n)$. But, $\si^\om(a_n)$ commutes
with $(1 \ot p_L)$ as well. Because $\| \si^\om_t(a_n) - a_n \| \recht
0$ uniformly on compact subsets of $\R$, it follows that $\|[(\io \ot
\mu)(x_L),a_n] \| \recht 0$. The same holds for $b_n$ and $c_n$. From
Lemma \ref{lemma.barnett}, we get that $(\io \ot \mu)(x_L) \in \C$. This holds
for all $\mu \in \B(L^2(\R))_*$ and we conclude that $x_L \in 1 \ot
\B(L^2(\R))$ for all compact subsets $L \subset \R$. We finally get
that $x \in 1 \ot \B(L^2(\R))$ as well. Because $x \in C$, it follows
that $x=1 \ot y$ with $y \in \rho(\R)\dpr$. Since, with the notations of the
proof of item \ref{e.a}), $\{(\mu \ot \io)\si^\om(N_1) \mid \mu \in
N_* \}\dpr = L^\infty(\R)$, we get $y \in \rho(\R)\dpr \cap
L^\infty(\R)' = \C$. This proves our claim.

\ref{e.c}) Suppose that $t_0 \in T(M,\de)$. Exactly as in item
\ref{e.b}), we can start with the strictly outer action $\al$ of
$(M,\de)$ on $N$ obtained in item \ref{e.a}) and next, extend this
action, to an action $\ga$ of $(M,\de)$ on $$C_{t_0}:= N
\rtimes_{(\si^\om_{nt_0})} \Z \; .$$
We can consider $C_{t_0}$ as an intermediate subalgebra: $\si^\om(N)
\subset C_{t_0} \subset C$, where $C$ denotes the continuous core of $N$. Because
we have seen already that $C \cap \si^\om(N^\al)' = \C$, it follows
that $C_{t_0} \cap \si^\om(N^\al)' = \C$. Hence, $\ga$ will be a
strictly outer action. To conclude the proof of item \ref{e.c}), we
have to show that $C_{t_0} \cong T_\sla$, with
$|t_0|=\frac{2\pi}{|\log \sla|}$.

As mentioned above, Shlyakhtenko has shown in Theorem 5.2 of
\cite{shl3} that the continuous core of $N:=\free{L^2(\R,K_\R), \sla_t}$
is isomorphic with $\cL(\F_\infty) \ot \B(K)$. More precisely, he has
given a model for the continuous core of $N$
that corresponds to the model for $\cL(\F_\infty) \ot \B(K)$ given by
R\u adulescu \cite{rad1} and such that the dual action $(\te_s)$ on the continuous
core is precisely the trace-scaling action of R\u adulescu. Moreover,
as explained just before Proposition 6.9 in \cite{shl} (using the
result of R\u adulescu proven in \cite{rad2}), the discrete
core of $T_\sla$ with its dual action $(\mu_n)$ of $\Z$ can be
described by the same model of $\cL(\F_\infty) \ot \B(K)$ and $\mu_n$
corresponds to $\te_{n \log \sla}$. Because $T_\sla$ is isomorphic
with the fixed point algebra of the dual action on the discrete core
of $T_\sla$, we get
$$T_\sla \cong \{ x \in N \rtimes_{(\si_t^\om)} \R \mid \te_{\log
  \sla}(x) = x \} \; .$$
It is clear that this last algebra is isomorphic with our $C_{t_0}$.

\ref{e.d}) From \ref{cor.faithaction}, we get a faithful action of
$(M,\de)$ on $\cL(\F_\infty)$ leaving invariant the trace. Taking the free product with a trivial
action on $\cL(\F_2)$, we immediately get a strictly outer action of
$(M,\de)$ on $\cL(\F_\infty)$ leaving invariant the trace.

\ref{e.e}) We start once again with the strictly outer action $\al$ of
$(M,\de)$ on $N$ obtained in item \ref{e.a}). As in the proof of item
\ref{e.b}), we can extend $\al$ to a strictly outer action $\ga$ of
$(M,\de)$ on $N \rtimes_{(\si_t^\om)} \Gamma$. Because $N$ is a type III$_1$
factor (\cite{shl}, Theorem 6.10), $N \rtimes_{(\si_t^\om)} \Gamma$ is a type
III$_0$ factor whose flow of weights is the action of $\R \subset \Gammah$
on the compact group $\Gammah$ by translations.
\end{proof}

We recall that a l.c.\ quantum group is said to be \emph{unimodular}
if the left invariant weight is also right invariant.

\begin{corollary}
If the dual of a l.c.\ quantum group $(M,\de)$ is unimodular, then
$(M,\de)$ can act strictly outerly on $\cL(\F_\infty) \ot \B(K)$. This holds in
particular for discrete quantum groups.

If a compact quantum group acts strictly outerly on a II$_\infty$ factor, then $\tau_t=\io$ for all $t \in \R$.
\end{corollary}
\begin{proof}
Let $(M,\de)$ be a l.c.\ quantum group. If the dual of $(M,\de)$ is
unimodular, it follows from Proposition 2.4 in \cite{VVD} that $\tau_t
= \Ad \sde^{-it/2}$, where $\sde$ is the modular element of $(M,\de)$, i.e.\ the Radon-Nikodym derivative between the left and right invariant weight.
Hence, $(M,\de)$ can act on $\cL(\F_\infty) \ot \B(K)$, because $\de(\sde^{it}) = \sde^{it} \ot \sde^{it}$. It is well known that a compact quantum group
is unimodular and so, our result can be applied to discrete quantum groups.

If a compact quantum group acts strictly outerly on a II$_\infty$ factor, we get a strictly positive, self-adjoint operator $\rho$ affiliated
with $M$ such that $\de(\rho^{it}) = \rho^{it} \ot \rho^{it}$ for all $t \in \R$ and $\tau_t = \Ad \rho^{it}$. Because $\de(\rho) = \rho \ot \rho$,
it follows with the same proof as in Proposition 7.10 of \cite{KV1} that $\rho$ is affiliated (in the C$^*$-algebraic sense) with the C$^*$-algebra
of $(M,\de)$. This C$^*$-algebra is unital and it follows that $\rho$ is bounded. Because $\de(\rho)=\rho \ot \rho$, we get $\rho = 1$. So, $\tau_t =
\io$ for all $t \in \R$.
\end{proof}

We give some more information on the $T$ invariant $T(M,\de)$ and calculate it for the bicrossed product quantum groups of \cite{VV,BSV}.

The following result is not true for the usual $T$ invariant of von Neumann algebras.
\begin{proposition}
Let $(M,\de)$ be a l.c.\ quantum group (and we still assume that $M_*$ is separable). Then $T(M,\de)$ is either a countable subgroup of $\R$, either the whole of
$\R$. In the latter case, there exists a one-parameter family $\rho^{it}$ of unitaries in $M$ satisfying $\de(\rho^{it}) = \rho^{it} \ot \rho^{it}$
and $\tau_t = \Ad \rho^{it}$ for all $t \in \R$.
\end{proposition}
\begin{proof}
We define the intrinsic group of $(M,\de)$.
$$G_M := \{ u \in M \mid u \quad\text{is unitary and}\quad \de(u) = u \ot u \} \; .$$
Equipped with the strong$^*$ topology, $G_M$ is a l.c.\ group. We know that for $u \in G_M$, $\tau_s(u)=u$ for all $s \in \R$. So, if $v \in G_M$ and
$\Ad v = \tau_s$ for some $s \in \R$, then $v \in \cZ(G_M)$, the center of $G_M$. We define the l.c.\ abelian group $K$ as
$$K:= \{ (u,t) \in \cZ(G_M) \times \R \mid \Ad u = \tau_t \}$$
and the continuous homomorphism $\eta : K \recht \R : \eta(u,t) = t$. By definition and the remark above, $T(M,\de) = \eta(K)$.

Define $K_0 = \operatorname{Ker} \eta$ and consider the dual,
continuous homomorphism $\hat{\eta} : \R \recht \widehat{K/K_0}$ with
dense range. So, $\widehat{K/K_0}$ is connected and hence, isomorphic
to $\R^n \times L$ for $n \in \N \cup \{0\}$ and $L$ a compact,
abelian group. We conclude that either $\widehat{K/K_0}$ is $\R$ or a
compact group. Then, $K/K_0$ is either $\R$ or a discrete group that
is second countable by construction. So, we have proved the first part of the proposition.

If $T(M,\de) = \R$, we get that $\eta$ is a continuous, surjective
homomorphism from $K$ onto $\R$. From \cite{HR}, Theorem 24.30, we know that $K$ is
isomorphic with $\R^n \times L$, where $n \in \N \cup \{0\}$ and $L$ is a l.c.\ abelian group
containing a compact open subgroup. Suppose that the restriction of
$\eta$ to $\R^n$ is not surjective. Then, $\eta(\R^n) = \{0\}$ and the
restriction of $\eta$ to $L$ is surjective. The image of the compact
open subgroup of $L$ is a compact subgroup of $\R$ and hence,
trivial. So, the image of $L$ is at most countable, contradicting the
surjectivity of $\eta$. Hence, the restriction of $\eta$ to $\R^n$ is
surjective and we find a continuous homomorphism $\mu : \R \recht K$
such that $\eta \mu = \io$. So, we are done.
\end{proof}

We now compute $T(M,\de)$ for bicrossed product l.c.\ quantum groups. We use the conventions of \cite{BSV}. So, we are given a (second countable)
l.c.\ group $G$ with two closed subgroups $G_1,G_2$ such that $G_1 \cap G_2 = \{e\}$ and $G_1 G_2$ has a complement of Haar measure zero in $G$. We
say that the pair $G_1,G_2$ is a matched pair. We define for almost all $x \in G$, $p_i(x) \in G_i$ such that $x=p_1(x)p_2(x)$. In this way, we identify (as
measure spaces) $G_1$ and $G/G_2$. We define the von Neumann algebra $M = G_2 \ltimes L^\infty(G/G_2)$. From \cite{VV}, we know that $M$ is the
underlying von Neumann algebra of a l.c.\ quantum group $(M,\de)$, the bicrossed product of $G_1$ and $G_2$. The underlying von Neumann algebra of
the dual $(\Mh,\deh)$ is given by $L^\infty(G_1 \backslash G) \rtimes G_1$. We have the following result.

\begin{proposition} \label{prop.Tbicross}
If $(M,\de)$ is the above bicrossed product l.c.\ quantum group, then
$$T(M,\de) = \{ t \in \R \mid \; \text{The character}\; \sde_2^{it} \;\text{on}\; G_2 \;\text{can be extended to a continuous character of}\; G
\} \; ,$$ where $\sde_2$ is the modular function of the l.c.\ group $G_2$.
\end{proposition}
\begin{proof}
We represent $M$ as $\bigl(\cL(G_2) \ot 1 \;\cup\; \al(L^\infty(G_1)) \bigr)^{\prime\prime}$ on $L^2(G_2 \times G_1)$, where $\cL(G_2)$ is the group
von Neumann algebra generated by the left translations and $\al(F)$ is the multiplication operator with the function $\al(F)(s,g) = F(p_1(sg))$ for
$g \in G_1, s \in G_2$, whenever $F \in L^\infty(G_1)$. From Proposition 4.16 in \cite{VV}, we know that $\tau_t$ is implemented by the
multiplication operator
$$P^{it}(s,g) = \sde^{it}(g^{-1} p_1(sg) ) \; \sde^{it}_1(g p_1(sg)^{-1}) \; \sde^{it}_2(s^{-1} p_2(sg)) \; ,$$
where $\sde_1,\sde_2,\sde$ denote the modular functions of $G_1,G_2,G$. Given our representation of $M$, it is clear that $1 \ot L^\infty(G_1)
\subset M'$. Hence, $\tau_t$ is implemented by the multiplication
operator with the function
$$(s,g) \mapsto \al(\sde^{it} \; \sde_1^{-it})(s,g) \;
\sde_2^{it}(s^{-1} p_2(sg)) \; .$$
We also know from Proposition 4.1 in \cite{BV} that the group-like unitaries of $(M,\de)$ are precisely of the form $\al(\om) (\sla_s
\ot 1)$, where $\om$ is a character on $G_1$ and $s \in G_2$ is such that $p_2(sg)=s$ for all $g \in G_1$.
Then, it is not difficult to conclude that $t \in T(M,\de)$ if and only if there exists a character $\om$ on $G_1$ such that the multiplication
operator
$$K(s,g) = \om(p_1(sg)) \; \sde_2^{it}(s^{-1} p_2(sg)) $$
commutes with $\cL(G_2) \ot 1$. This holds if and only if
$$\om(p_1(sg)) \; \sde_2^{it}(p_2(sg)) = \om(g) \; \sde_2^{it}(s) \quad\text{almost everywhere} \; ,$$
i.e.\ if and only if there exists a (measurable, hence continuous) character $\mu$ on $G$ such that $\mu(gs) = \om(g) \; \sde_2^{it}(s)$.
So, we are done.
\end{proof}

In Example 5.4 of \cite{VV}, there is given a matched pair $G_1,G_2 \subset G$, such that $G \cong \operatorname{PSL}_2(\R)$ and $G_1$ is the
$(ax+b)$-group. Interchanging $G_1$ and $G_2$, the bicrossed product $(M,\de)$ satisfies $T(M,\de)=\{0\}$.

We also remark that, for an arbitrary l.c.\ quantum group $(M,\de)$, the scaling group $(\tau_t)$ satisfies $\de \tau_t=(\tau_t \ot \tau_t)\de$.
Hence, $\R$ acts by quantum group automorphisms on $(M,\de)$. If $\Gamma$ is any subgroup of $\R$, that we equip with the discrete topology, the
crossed product $\Gamma \kruisje{(\tau_t)} M$ carries a natural l.c.\ quantum group structure: the comultiplication on $M$ is extended by putting
$\de(\sla_x) = \sla_x \ot \sla_x$. It is now obvious, that if $T(M,\de)=\{0\}$, then $T(\Gamma \ltimes M,\de) = \Gamma$. So, we can obtain any
subgroup of $\R$ as the $T$ invariant of a l.c.\ quantum group and any countable subgroup as the $T$ invariant of a l.c.\ quantum group with separable predual.

Another case where the invariant $T(M,\de)$ can be easily calculated,
is for certain compact quantum groups. We state the result for
Woronowicz' quantum $SU_q(n)$ group.

\begin{proposition}
Consider the compact quantum group $SU_q(n)$, $0 < q < 1$ and $n \geq
2$ \cite{Wor3}. Then,
$$T(SU_q(n)) = \Z \frac{\pi}{\log q} \; .$$
In particular, $SU_q(n)$ acts strictly outerly on the free Araki-Woods
factor $T_{q^2}$ of type III$_{q^2}$.
\end{proposition}
\begin{proof}
Write $(M,\de)=SU_q(n)$.
Since a group-like unitary $u \in M$, is a one-dimensional
corepresentation and since we know from \cite{Wor3} the
corepresentations of $(M,\de)$, we conclude that $(M,\de)$ has no
non-trivial group-like unitaries. Hence, $T(M,\de)$ consists of those
$t \in \R$ with $\tau_t = \io$. But, this is the same set as the set
of $t \in \R$ with $\tauh_t = \io$. We know that the scaling group on
the dual is implemented by the $F$-matrices. Hence, $\tauh_t = \io$ if
and only if $q^{2it}=1$.
\end{proof}

In the course of the proof of the Theorem \ref{thm.construction}, we needed the extension of an action $\al$ of a l.c.\ quantum group $(M,\de)$ on a von
Neumann algebra $N$ to an action on the crossed product $N \rtimes_{(\si^\te_t)} \R$, where $\te$ is an n.s.f.\ weight on $N$. Such extensions were
considered by Yamanouchi in \cite{Y} under the assumption that $\nabh = \sde^{-1}$, which implies in particular that $\tau_t = \Ad \sde^{-it}$, i.e.\
the scaling group is implemented by the modular element, which is group-like.

It is more natural to consider these extensions in the more general
setting where the scaling group is implemented by an arbitrary
group-like operator. In fact, we show that one cannot do better.

\begin{proposition}
Let $(M,\de)$ be a l.c.\ quantum group. The following are equivalent.
\begin{itemize}
\item Every action $\al : N \recht M \ot N$ of $(M,\de)$ on a von
  Neumann algebra $N$ can be extended to an action of $(M,\de)$ on
  the crossed product $N \rtimes_{(\si^\te_t)} \R$, whenever $\te$ is
  an n.s.f.\ weight on $N$.
\item There exists a strictly positive, self-adjoint operator $\rho$
  affiliated with $M$ and satisfying $\de(\rho^{it})=\rho^{it} \ot
  \rho^{it}$, $\tau_t = \Ad \rho^{it}$.
\end{itemize}
If the second condition is fulfilled, there is a canonical extension,
given by
$$\be \site = (\io \ot \site)\al \quad\text{and}\quad \be(1 \ot
\rho_t) = (\io \ot \site)(D_t) \; (\rho^{is} \ot 1 \ot \rho_s) \; ,$$
where
\begin{align*}
& \site : N \recht N \ot L^\infty(R) : \site(x)(t)=\site_t(x) \; ,
\quad (\rho_t \xi)(s) = \xi(s+t) \quad\text{for}\quad \xi \in L^2(\R)
\quad\text{and} \\ & N \rtimes_{(\si^\te_t)} \R = (\site(N) \cup 1 \ot
\rho_t)\dpr \subset N \ot \B(L^2(\R)) \; , \\ & D_t = [D \te \circ \al
: D \te]_t \in M \ot N \quad\text{is the Radon-Nikodym derivative.}
\end{align*}
\end{proposition}

\begin{proof}
Suppose first that every action $\al$ can be extended to an action on
the core. Consider the strictly outer action $\al$ of $(M,\de)$ on
$(N,\om) = (N_0,\om_0) \star (N_1,\om_1)$ constructed in the proof of
Theorem \ref{thm.construction}.\ref{e.a}), using the notation
introduced there. By assumption, there exists an action $\be : C
\recht M \ot C$ extending $\al$, where $C = N \rtimes_{(\si_t^\om)}
\R$. Since $\be$ extends $\al$, $\be$ is a faithful action. Moreover,
from the proof of Theorem \ref{thm.construction}.\ref{e.b}), we get
that $C \cap \si^\om(N^\al)' = \C$. Hence, $C \cap (C^\be)' = \C$. So,
$\be$ is strictly outer. Because $C$ is a II$_\infty$ factor, the
second condition of the proposition follows from Theorem
\ref{thm.restrictions}.\ref{en.b}).

Suppose, conversely, that the second condition is fulfilled. Let $\al
: N \recht M \ot N$ be an action of $(M,\de)$ on a von Neumann algebra
$N$ with n.s.f.\ weight $\te$. Define $\cU \in M \ot N \ot
L^\infty(\R)$ by $\cU(t) = (\rho^{-it} \ot 1)D_t^*$. If we define, for
$z \in N \rtimes_{(\si^\te_t)} \R \subset N \ot \B(L^2(\R))$, $\be(z) =
\cU(\al \ot \io)(z) \cU^*$. One verifies immediately that
$$\be \si^\te = (\io \ot \si^\te)\al \quad\text{and}\quad \be(1 \ot
\rho_t) = (\io \ot \si^\te)(D_t)(\rho^{it} \ot 1 \ot \rho_t) \; .$$
The cocycle property of $D_t$ implies immediately that $\be$ is an action.
So, we are done.
\end{proof}

\section{A technical lemma}

We now include two lemmas that were needed to establish the trivial
relative commutant properties $N \cap (N^\al)' = \C$ and $C \cap
\si^\om(N^\al)' = \C$ in the proof of Theorem
\ref{thm.construction}. Lemma \ref{lemma.barnett} below is a generalization of
Theorem 11 in \cite{Bar}. We do not need the full strength of the
lemma, but it could be interesting to state its most
general form.

The beautiful idea of the proof of Lemma \ref{lemma.barnett} is due to
G. Skandalis.

In Theorem 11 in \cite{Bar}, one considers a free product $(N,\om) = (N_1,\om_1)
\star (N_2,\om_2)$ (see \cite{voi}), a unitary $a \in N_1$ in the centralizer of
$\om_1$ satisfying $\om_1(a) = 0$ and unitaries $b,c \in N_2$ in the
centralizer of $\om_2$ satisfying $\om_2(b)=\om_2(c)=\om_2(cb^*)=0$, to obtain the inequality
$$\|x-\om(x) 1\|_2 \leq 14 \max \{ \|[x,a]\|_2 , \|[x,b]\|_2 ,
\|[x,c]\|_2 \} \; , \; \; \text{for all}\;\; x \in N \; .$$
Here, and below, we use the $L^2$-norm $\|x\|_2 = \om(x^*x)$.

We will prove a same kind of inequality for general elements $a,b,c$
keeping track of their non-unitarity and non-invariance under the
modular group. This will allow to conclude the relative
commutant properties needed above (in cases where the centralizer is trivial).

\begin{lemma} \label{lemma.barnett}
Let $N_{1,2}$ be von Neumann algebras with faithful normal states $\om_{1,2}$.
Let $(N,\om) = (N_1,\om_1) \star (N_2,\om_2)$. Let $a \in
N_1$ and $b,c \in N_2$. Suppose that $a,b,c \in
D(\si^\om_{i/2})$. Then, for all $x \in N$,
$$\|x - \om(x)1 \|_2 \leq D(a,b,c) \, \max \{ \|[x,a]\|_2 , \|[x,b]\|_2 ,
\|[x,c]\|_2 \} + E(a,b,c) \, \|x\|_2 \; ,$$
where $D(a,b,c) = 3\|a\|^3 + 2 \|b\|^3 + 2 \|c\|^3$ and $E(a,b,c) = 3
C(a) + 2 C(b) + 2 C(c) + 6 |\om(cb^*)| \, \|cb^*\|$ with
$$C(a) = \|a\|^3 \, \|\si^\om_{i/2}(a) - a \| + \|a\|^2 \, \|a^*a-1\|
+ 2(1+\|a\|^2)\,\|aa^*-1\| + 3 |\om(a)| \, \|a\| \; .$$
Hence, if there exists a sequence $a_n$
in $N_1$ and sequences $b_n,c_n$ in $N_2$ such that
$a_n,b_n,c_n$ satisfy
\begin{align*}
& \|a_n^* a_n - 1 \| \recht 0 \; , \quad \|a_n a_n^* - 1 \| \recht 0
\; , \\ & a_n \in D(\si^\om_{i/2}) \quad\text{and}\quad \|
\si^\om_{i/2}(a_n) - a_n \| \recht 0 \; , \\ & \om(a_n) \recht 0 \;
.
\end{align*}
and such that $\om(c_n b_n^*) \recht 0$, then
$$N \cap \{ a_n,b_n,c_n \mid n \in \N \}' = \C \; .$$
Moreover, every central sequence in $N$ is trivial.
\end{lemma}
\begin{proof}
Take $a \in N_1$ and $b,c \in N_2$ such that $a,b,c \in
D(\si^\om_{i/2})$. Choose $x \in N$. Replacing $x$ by $x - \om(x)1$
(which reduces the $L^2$-norm), we may suppose that $\om(x) = 0$.
Represent $N_1$, resp.\ $N_2$, on its GNS-space
$(H_1,\xi_1)$, resp.\ $(H_2,\xi_2)$. Let $(H,\xi) = (H_1,\xi_1) \star
(H_2,\xi_2)$.
Recall that
$$H= \C \xi \oplus (\Honebol \ot H(2,l)) \oplus (\Htwobol \ot H(1,l))
\; ,$$
with
\begin{align*}
& H(2,l) = \C \xi \oplus \Htwobol \oplus (\Htwobol \ot \Honebol) \oplus
(\Htwobol \ot \Honebol \ot \Htwobol) \oplus \cdots
\quad\text{and}\\ & H(1,l) = \C \xi \oplus \Honebol \oplus (\Honebol \ot \Htwobol)
\oplus (\Honebol \ot \Htwobol \ot \Honebol) \oplus \cdots \; ,
\end{align*}
where $\Hbol_i = H_i \ominus \C \xi_i$.

Define $\eta:=x \xi$.
Because $\langle \eta,\xi \rangle = 0$, we write $\eta = \mu +
\ga$ with $\mu \in \Honebol \ot H(2,l)$ and $\ga \in \Htwobol
\ot H(1,l)$. Define, for $\zeta \in H$ and $y \in N$, $\zeta \cdot y
:= Jy^*J \zeta$ and observe that $(z \xi) \cdot \si^\om_{i/2}(y) = zy
\xi$ for $y \in D(\si^\om_{i/2})$.

Define $\etatil = a^* \cdot \eta \cdot a$. We observe that
$$\|\etatil - \eta \| \leq
\|a\| \, \|\si^\om_{i/2}(a) - a \| \, \|x\|_2 +
\|a\| \, \|[x,a]\|_2 + \|a^*a-1\| \, \|x\|_2 \; .$$
Writing $\mutil = a^* \cdot \mu \cdot a$ and $\gatil = a^* \cdot \ga
\cdot a$, we have $\etatil = \mutil + \gatil$. It is obvious that
$$|\langle \mutil,\gatil \rangle| \leq (1+\|a\|^2) \, \|aa^*-1\| \,
\|x\|_2^2 \quad\text{and}\quad \bigl| \, \|\gatil\|^2 - \|\ga\|^2 \,
\bigr| \leq  (1+\|a\|^2) \, \|aa^*-1\| \,
\|x\|_2^2 \; .$$
If $P_2$ denotes the projection onto $\Htwobol \ot H(1,l)$ and because
$a \in N_1$, we know
that $$\|P_2 a^* Ja^*J \ga\| \leq |\om(a)| \, \|Ja^*J \ga\| \leq
|\om(a)| \, \|a\| \, \|x\|_2 \; .$$
So, we conclude that $|\langle \gatil,\ga \rangle| \leq |\om(a)| \,
\|a\| \, \|x\|_2^2$.

Write $\zetatil = \eta - \ga - \gatil$. With the inequalities obtained
so far, we can estimate all the inner products between $\zetatil,\ga$
and $\gatil$. Since we also estimate the difference between
$\|\gatil\|^2$ and $\|\ga\|^2$, we arrive at
$$\|\eta\|^2 = \|\zetatil + \ga + \gatil\|^2 \geq
\|\zetatil\|^2 +2 \|\ga\|^2 - \|a\|^3 \, \|x\|_2 \, \|[x,a]\|_2 - C(a)
\|x\|_2^2 \; ,$$
where $C(a)$ is as in the statement of the lemma. Since $\|\eta\|^2 =
\|\mu\|^2+\|\ga\|^2$, we arrive at
\begin{equation} \label{eq.firstcrucial}
\|\mu\|^2 \geq \|\ga\|^2 - \|a\|^3 \, \|x\|_2 \, \|[x,a]\|_2 - C(a)
\|x\|_2^2 \; .
\end{equation}

Next, we use quite analogously $b$ and $c$. Define $\eta'=b^* \cdot
\eta \cdot b$, $\eta\dpr = c^* \cdot \eta \cdot c$ and analogously
$\mu',\mu\dpr,\ga'$ and $\ga\dpr$. Write $\zeta' = \eta - \mu - \mu' -
\mu\dpr$. We can estimate in exactly the same way as above, all the
inner products between $\zeta',\mu,\mu'$ and $\mu\dpr$. It will not be
surprising that we arrive in this way at
\begin{equation} \label{eq.secondcrucial}
\|\ga\|^2 \geq 2 \|\mu\|^2 - \|b\|^3 \, \|x\|_2 \, \|[x,b]\|_2
- \|c\|^3 \, \|x\|_2 \, \|[x,c]\|_2 - (C(b)+C(c)+3 |\om(cb^*)|\,
\|cb^*\|) \|x\|_2^2 \; .
\end{equation}
Combining Inequalities \eqref{eq.firstcrucial} and
\eqref{eq.secondcrucial} and using that $\|x\|_2 =
(\|\mu\|^2+\|\ga\|^2)/\|x\|_2$, we precisely arrive at the inequality
stated in the lemma.

Suppose now that we have sequences $a_n,b_n,c_n$ as stated in the
lemma. It is clear that $D(a_n,b_n,c_n)$ will remain bounded, while
$E(a_n,b_n,c_n)$ converges to zero. So, we have a constant $D \geq 0$
and a sequence $\kappa_n \geq 0$ converging to zero, such that
$$\|x - \om(x)1 \|_2 \leq D \max \{ \|[x,a_n]\|_2 , \|[x,b_n]\|_2 ,
\|[x,c_n]\|_2 \} + \kappa_n \, \|x\|_2 \; ,$$
for all $x \in N$ and $n \in \N$.

If $x \in N \cap \{ a_n,b_n,c_n \mid n \in \N \}'$, we immediately get
that $x=\om(x)1 \in \C$.

Let $x_n$ be a central sequence in $N$ (i.e.\ a bounded sequence in $N$
such that $\| [x_n,a]\|_2 , \| [x^*_n,a]\|_2 \recht 0$ for all $a \in
N$). Choose $\vareps > 0$. Take $n$ such that $\kappa_n \, \|x_m\|_2 <
\vareps/2$ for all $m$. Next, take $m_0$ such that for all $m \geq m_0$,
$D \max \{ \|[x_m,a_n]\|_2 , \|[x_m,b_n]\|_2 ,
\|[x_m,c_n]\|_2 \} < \vareps/2$. It follows that $\|x_m - \om(x_m)1
\|_2 < \vareps$ for all $m \geq m_0$. We can do the same thing with
$x_m^*$ and conclude that the sequence $x_n$ is trivial.
\end{proof}

The following lemma is very easy to prove.

\begin{lemma} \label{lemma.free}
The free Araki-Woods factor $N=\free{L^2(\R,\R),\sla_t}$ with free
quasi-free state $\om$ contains a
sequence of unitaries $u_n$ such that $u_n$ is analytic w.r.t.\
$(\si^\om_t)$, $\|\si^\om_z(u_n) - u_n \| \recht 0$ uniformly on
compact subsets of $\C$ and $\om(u_n) \recht 0$.
\end{lemma}
\begin{proof}
By Fourier transformation, we consider rather $K_\R = \{ \xi \in
L^2(\R) \mid \xi(-x) = \overline{\xi(x)} \}$, with orthogonal
transformations $(U_t \xi)(x) = \exp(itx) \xi(x)$. This corresponds to
the involution $(T \xi)(x) = \exp(-x/2) \overline{\xi(-x)}$ on $L^2(\R)$.

Take unit vectors
$\mu_n$ in $K_\R$, bounded and with support in $[-1/n,1/n]$. Recall that generators
for the free Araki-Woods factor $N$ can be written as $s(\xi) =
(\ell(\xi) + \ell(T \xi)^*)/2$. Also, if $\xi$ is a bounded, compactly
supported function in $L^2(\R)$, then $s(\xi)$ is analytic w.r.t.\
$(\si^\om_t)$ and $\si^\om_z(s(\xi)) = s(\xi_z)$, where $\xi_z(x) =
\exp(izx) \xi(x)$.

Put $x_n:= (s(\mu_n) + s(\mu_n)^*)/2
\in N$. Put also $y_n = (\ell(\mu_n)+\ell(\mu_n)^*)/2 \in
\B(\cF(L^2(\R)))$. We know that $y_n$ has the semi-circular
distribution with respect to the vacuum state. Take a real number $k$
such that $\int_{-1}^1 \exp(ikt) \sqrt{1-t^2} \, dt = 0$. Put $u_n =
\exp(ikx_n)$. It is clear that $u_n$ are unitaries in $N$, analytic w.r.t.
$(\si^\om_t)$ and $\| \si^\om_z(u_n) - u_n \| \recht 0$ uniformly on
compact subsets of $\C$. Since $\|x_n - y_n\| \recht 0$, we get
$\|u_n - \exp(iky_n)\| \recht 0$, which yields $\om(u_n) \recht 0$ by
our choice of $k \in \R$.
\end{proof}

\section{Strictly outer actions of locally compact groups on \\ the
  hyperfinite II$_1$ factor}

In this section, we construct strictly outer actions of ordinary l.c.\
groups on the hyperfinite II$_1$ factor. We first provide a general
construction procedure, based on a general strict outerness
result. Next, we give a more geometrical construction for linear groups.

\begin{theorem} \label{thm.outerR}
Let $N$ be a factor with a faithful state $\om$ and suppose that a
l.c.\ group $G$ acts faithfully on $N$ by automorphisms $(\al_g)$
leaving invariant $\om$, i.e.\ such that
$\al_g = \io$ implies $g=e$ and $\om \al_g =\om$.

Denote $(N_\infty,\om_\infty) = \bigotimes_{n=1}^\infty (N,\om)$, which is a factor
equipped with the faithful state $\om_\infty$. Let $G$ act diagonally
on $N_\infty$.

Then, the action of $G$ on $N_\infty$ is strictly outer.
\end{theorem}

\begin{proof}
Write $N_k = \bigotimes_{n=1}^k (N,\om)$.

For every $n \in \N$, it is easy to define an automorphism $i_n$ of $N_\infty$
such that
$$i_n(x_1 \ot \cdots \ot x_n \ot \cdots) = x_n \ot x_1 \ot \cdots \ot
x_{n-1} \ot x_{n+1} \ot \cdots \; .$$
One defines $i_n$ on all $N_k$, $k \geq n$, and one can extend to
$N_\infty$ because $i_n$ preserves $\om_\infty$.

We have an obvious isomorphism $\Psi: N \ot N_\infty \recht N_\infty$
given by $\Psi(x_0 \ot (x_1 \ot x_2 \ot \cdots)) = x_0 \ot x_1 \ot
\cdots$. We write $j(z) = \Psi(1 \ot z)$ and claim that $i_n(z) \recht
j(z)$ strongly$^*$ for all $z \in N_\infty$.

To prove our claim, let $(H_\infty,\xi_\infty) = \bigotimes_{n=1}^\infty (H,\xi_\om)$ be a
GNS-construction for $(N_\infty,\om_\infty)$ with cyclic and
separating vector $\xi_\infty$. We define unitaries $U_n$ on
$H_\infty$ by the formula $U_n x \xi_\infty = i_n(x) \xi_\infty$. We
also define an isometry $V$ on $H_\infty$ given by $V x \xi_\infty =
j(x) \xi_\infty$. By definition, it is clear that $U_n x \xi_\infty
\recht V x \xi_\infty$ whenever $x \in N_k$ for some $k \in
\N$. Because the sequence $U_n$ is bounded, it follows that $U_n
\recht V$ strongly. Hence, for $z \in N_\infty$ and $y \in N_\infty'$,
we get
$$i_n(z) y \xi_\infty = y i_n(z) \xi_\infty = y U_n z \xi_\infty \recht y V z \xi_\infty = y
j(z) \xi_\infty = 
j(z) y \xi_\infty \; .$$
Because the sequence $i_n(z)$ is bounded, we conclude that $i_n(z)
\recht j(z)$ strongly for all $z \in N_\infty$. Because $i_n$ and $j$
are $^*$-homomorphisms, the convergence is strong$^*$. This proves our
claim.

We realize $G \ltimes N_\infty = (\al(N_\infty) \cup \cL(G) \ot
1)\dpr \subset \B(L^2(G)) \ot N_\infty$, where $\al : N_\infty \recht
L^\infty(G) \ot N_\infty : \al(x)(g) = \al_{g^{-1}}(x)$.

Let now $a \in G \ltimes N_\infty \cap \al(N_\infty)'$.
Using the realization of $G \ltimes N_\infty \subset \B(L^2(G)) \ot N_\infty$ and the fact that $i_n$ and $j$ commute with
the action of $G$, we find that
$$(\io \ot i_n)(a) \in G \ltimes N_\infty \cap \al(N_\infty)'$$
for all $n$ and hence,
$$(\io \ot j)(a) \in G \ltimes N_\infty \cap \al(N_\infty)' \; .$$
Applying the isomorphism $\io \ot \Psi^{-1}$, we conclude that
$a_{13}$ and $\al(x) \ot 1$ commute in $\B(L^2(G)) \ot N \ot N_\infty$
for all $x \in N$. So, $a$ and $(\io \ot \mu)\al(x) \ot 1$ commute in 
$\B(L^2(G)) \ot N_\infty$ for all $x \in N, \mu \in N_*$. The
self-adjoint family of functions
$$\{ g \mapsto \mu(\al_{g^{-1}}(x)) \mid  \mu \in N_*, x \in N \}$$
separates the points of $G$, which yields
$$\{ g \mapsto \mu(\al_{g^{-1}}(x)) \mid  \mu \in N_*, x \in N \}' =
L^\infty(G) \; .$$
It follows that $a \in L^\infty(G) \ot N_\infty$. From the duality for
crossed products, we get that
$$a \in G \ltimes N_\infty \cap L^\infty(G) \ot N_\infty =
\al(N_\infty) \; .$$ But, $a \in \al(N_\infty)'$, which gives $a \in
\C 1$ and we are done.
\end{proof}

In \cite{Bla}, Blattner constructed a canonical action of any l.c.\
group on the hyperfinite II$_1$ factor. His construction is as
follows. Let $H_\R$ be a real Hilbert space and consider its Clifford
algebra $\Cl$: $\Cl$ is the $^*$-algebra generated by the self-adjoint elements $c(\xi)$, $\xi \in
H_\R$ with relations $c(\xi) c(\eta) + c(\eta) c(\xi) = 2 \langle
\xi,\eta \rangle$. For every orthogonal transformation $u$ on $H_\R$,
we have an automorphism $\al_u$ of $\Cl$ such that
$\al_u(c(\xi)) = c(u \xi)$. In particular, taking $u \xi = -\xi$ for
all $\xi \in H_\R$, we obtain a $\Z / 2\Z$-grading of $\Cl$. On $\Cl$
there exists a unique trace $\tau$ such that $\tau(1) = 1$ and
$\tau(a) = 0$ if $a$ has odd degree. Using a GNS-representation for
this trace, we define a von Neumann algebra, which happens to be the
hyperfinite II$_1$ factor $\cR$ if $H_\R$ is of separable infinite
dimension. So, whenever such a real Hilbert space $H_\R$ is fixed, we
consider $$\cR = \{ c(\xi) \mid \xi \in H_\R \}^{\prime\prime} \; .$$
It is clear that, for any orthogonal transformation $u$ on $H_\R$,
$\al_u$ extends to an automorphism of $\cR$, still denoted by $\al_u$.

Let now $G$ be a l.c.\ group and let $(u_g)$ be a continuous
representation of $G$ by orthogonal transformations of a real Hilbert
space $H_\R$ of infinite separable dimension. Defining $\al_g :=
\al_{u_g}$, we get an action of $G$ by automorphisms of $\cR$. If the
representation $(u_g)$ is faithful, the automorphism group $(\al_g)$
is clearly faithful as well.

Also, observe that taking a direct sum of a family of orthogonal
transformations yields an action of $\cR$, which is isomorphic to the
diagonal action on the tensor product of the family of copies of $\cR$
associated to the given family of real Hilbert spaces.

Hence, we obtain the following corollary of Theorem \ref{thm.outerR}.

\begin{corollary} \label{cor.groupouter}
Every l.c.\ group $G$ can act strictly outerly on the hyperfinite
II$_1$ factor. In particular, Blattner's action of $G$ on the
hyperfinite II$_1$ factor associated as above with an infinite direct
sum of copies of the regular representation, is a strictly outer action.
\end{corollary}

Finally, we explain how strictly outer actions of linear groups on the
hyperfinite II$_1$ factor can be obtained in an alternative, more intuitive and
geometric way.

Let $n \in \N$. Define $G = \SL(n+2,\C)$ and consider the subgroups
\begin{align*}
H &= \Biggl\{\begin{pmatrix} 1_{n,n} & 0_{n,2} \\ 0_{2,n} & \begin{matrix} 1 & k
    \\ 0 & 1 \end{matrix} \end{pmatrix}  \Bigg| \; k \in \Z \Biggr\}
\; , \quad  K = \Biggl\{\begin{pmatrix} \SL(n,\C) & 0_{n,2} \\
    0_{2,n} & \begin{matrix} 1 & z \\ 0 & 1 \end{matrix} \end{pmatrix}
\Bigg| \; z \in \C \Biggr\} \quad\text{and} \\
L &= \Biggl\{\begin{pmatrix} \SL(n,\C) & 0_{n,2} \\
    0_{2,n} & \begin{matrix} 1 & k \\ 0 & 1 \end{matrix} \end{pmatrix}
\Bigg| \; k \in \Z \Biggr\} \; .
\end{align*}
Let $\Ga$ be a lattice in $G$ (i.e.\ a discrete subgroup with
finite covolume). Write $X = G / \Gamma$, with its finite measure
invariant under the action of $G$.

We first claim that the action of $K$ on $X$ is (measure
theoretically) free, in the sense
that for almost all $x \in X$, the stabilizer $S_x \subset K$ is
trivial. So, we have to prove that
$$\bigcup_{\gamma \in \Ga, \gamma \neq e} \{ x \in G \mid x \gamma
x^{-1} \in K \}$$
has measure zero in $G$. Because $\Ga$ is countable, it is enough to
show that, for every $\gamma \in \Gamma \setminus \{e\}$, $\{ x \in G \mid x \gamma x^{-1} \in K \}$ has
measure zero.
This clearly is an algebraic subvariety of $G$, which is
not the whole of $G$, because the center of $G$ intersects $K$
trivially. Hence, it is a set of measure zero. We conclude a fortiori
that the action of $L$ on $X$ is free in the same measure theoretic sense.

In particular, the action of $H$ on $X$ is free and it is ergodic by
Moore's ergodicity theorem, see e.g.\ \cite{Zim}. So, we can define $\cR = H \ltimes
L^\infty(X)$, which is the hyperfinite II$_1$ factor by amenability of
$\Z$, ergodicity and freeness.

We next have a natural action of $\SL(n,\C)$ on $\cR$ such that
$$\SL(n,\C) \ltimes \cR = L \ltimes L^\infty(X) \; .$$
It follows that $\SL(n,\C) \ltimes \cR \cap \cR' \subset L \ltimes
L^\infty(X) \cap L^\infty(X)' = L^\infty(X)$ because of the freeness
of the action of $L$ on $X$, see \cite{Sau}. But then, $\SL(n,\C) \ltimes \cR \cap \cR'
\subset L^\infty(X) \cap \cR' = \C$ because of the ergodicity of the
action of $H$ on $X$.

Hence, we have found a strictly outer action of $\SL(n,\C)$ on the
hyperfinite II$_1$ factor $\cR$ for any $n \in \N$. Because the restriction of a strictly
outer action to a closed subgroup is a strictly outer action, we have
found a strictly outer action on $\cR$ of any linear group.

\section{Strictly outer actions of locally compact quantum groups on
  \\ injective factors}

First, we have shown that every l.c.\ quantum group can act strictly
outerly on a free Araki-Woods factor. These factors are far from being
injective. Secondly, we have seen that every l.c.\ group can act
strictly outerly on the hyperfinite II$_1$ factor. From Theorem
\ref{thm.restrictions}, we know that not all l.c.\ quantum groups can
act strictly outerly on a II$_1$ factor. Nevertheless, we ask the
natural question when it is possible to act strictly outerly on an
injective factor of arbitrary type.

In this elementary section, we prove a general result giving a
necessary condition for the possibility to act strictly outerly on an
injective factor. We prove that for bicrossed product quantum groups,
this condition is sufficient as well.

Given a corepresentation of a l.c.\ quantum group in a factor $N$, we
consider the infinite tensor product of the associated inner action on
$N$. We give a sufficient condition for its strict outerness that will
allow us, in the next two sections, to deal with strictly outer
actions of compact and discrete quantum groups on injective factors.

We finally show that he possibility
of acting strictly outerly on an injective factor is stable under
cocycle deformation of the quantum group.

\begin{definition}
Let $(M,\de)$ be a l.c.\ quantum group. A (not necessarily normal)
state $m : M \recht \C$ is called a left invariant mean, if $m \bigl(
(\om \ot \io)\de(x) \bigr) = \om(1) \, m(x)$ for all $x \in M$ and
$\om \in M_*$. We analogously define a right invariant mean. We call
$m$ an invariant mean, if $m$ is a left and right invariant mean.

A l.c.\ quantum group has a left invariant mean if and
only if it has an invariant mean, see e.g.\ Proposition 3 in \cite{dqv}.
\end{definition}

\begin{proposition} \label{prop.amenable}
Let $\al : N \recht M \ot N$ be a strictly outer action of a l.c.\
quantum group $(M,\de)$ on $N$. If the crossed product $M
\kruisje{\al} N$ is injective, $(M,\de)$ has an invariant mean. If $N$
is injective, the dual l.c.\ quantum group $(\Mh,\deh)$ has an
invariant mean.
\end{proposition}
\begin{proof}
It suffices to prove the first statement: the second follows by
considering the dual action on the crossed product and using the fact
that the double crossed product is $\B(H) \ot N$, which is injective
if $N$ is injective (see \cite{V} for the notion of dual action and
double crossed product).

So, suppose that $M \kruisje{\al} N$ is injective. Let $P : \B(H) \ot N
\recht M \kruisje{\al} N$ be a norm one projection. Take $z \in
M'$. Then, $P(z \ot 1) \in M \kruisje{\al} N$ and for all $x \in N$, we have
$$P(z \ot 1) \, \al(x) = P \bigl( (z \ot 1) \al(x) \bigr) = P\bigl(
\al(x) (z \ot 1) \bigr) = \al(x) \, P(z \ot 1) \; .$$
By strict outerness of $\al$, we get that $P(z \ot 1) \in \C$. Hence,
we can define a state $\mu$ on $M'$ such that $P(z \ot 1) = \mu(z) 1$
for all $z \in M'$.

Define $\rho : M' \recht M \ot M' : \rho(z) = W(1 \ot z)W^*$. Here,
$W$ is the multiplicative unitary associated with $(M,\de)$ and we
know that $(\Jh \ot J)W(\Jh \ot J) = W^*$ (see Preliminaries). Hence,
we get that $\rho(z) = (\Jh \ot J) \de(JzJ) (\Jh \ot J)$ for all $z
\in M'$. Let $Q$ be the norm one projection from $\B(H) \ot \B(H) \ot
N$ onto $\B(H) \ot (M \kruisje{\al} N)$ such that $(\om \ot \io \ot
\io)Q(z) = P((\om \ot \io \ot \io)(z))$ for all $z \in \B(H) \ot \B(H) \ot
N$ and $\om \in \B(H)_*$. We observe that, for $z \in M'$ and $\om \in \B(H)_*$,
\begin{align*}
\mu\bigl((\om \ot \io)\rho(z)\bigr) \, (1 \ot 1) &= P\bigl( (\om \ot
\io)\rho(z) \ot 1) = (\om \ot \io \ot \io)Q(\rho(z) \ot 1) \\ &=
(\om \ot \io \ot \io)Q(W_{12} (1 \ot z \ot 1) W_{12}^* ) \\ &=
(\om \ot \io \ot \io)\bigl( W_{12} \; Q(1 \ot z \ot 1) \; W_{12}^*
\bigr) = \om(1) \; \mu(z) \; (1 \ot 1) \; .
\end{align*}
We immediately conclude that $m(z) = \mu(Jz^*J)$ defines a left
invariant mean on $(M,\de)$.
\end{proof}

Combining several results, we obtain the following.

\begin{proposition}
Let $(M,\de)$ be a bicrossed product locally compact quantum group, with
$M = G_2 \ltimes L^\infty(G/G_2)$ as explained just before Proposition
\ref{prop.Tbicross}. Then, $(M,\de)$ can act strictly outerly on an
injective factor if and only if $(\Mh,\deh)$ has an invariant mean.
\end{proposition}
\begin{proof}
One implication follows from Proposition \ref{prop.amenable}. So,
suppose that $(\Mh,\deh)$ has an invariant mean. Since
$(L^\infty(G_2),\de_2)$ is a closed quantum subgroup of $(\Mh,\deh)$,
the restriction of the invariant mean on $(\Mh,\deh)$ gives an
invariant mean on $L^\infty(G_2)$. So, $G_2$ is an amenable l.c.\ group.
Using Corollary \ref{cor.groupouter}, we can take a strictly outer
action of $G$ on the hyperfinite II$_1$ factor $\cR$. From Proposition 6.1 in
\cite{BV}, we get that $(M,\de)$ can act strictly outerly on $G_2
\ltimes \cR$, where $G_2$ acts by the restriction of the action of $G$
on $\cR$. Since $G_2$ is amenable, $G_2 \ltimes \cR$ is injective.
\end{proof}

Taking into account Theorem \ref{thm.outerR} and since the tensor
product preserves injectiveness, it is a natural idea to
consider an infinite tensor product action to obtain strictly outer
actions on injective factors. But, due to the non-commutativity of the
algebra $M$, we cannot perform the tensor product of two arbitrary
actions. However, we can make the tensor product of two inner actions
because we can make the tensor product of corepresentations.

Let $(M,\de)$ be a l.c.\ quantum group and $U \in M \ot N$ a
corepresentation of $(M,\de)$ in the factor $N$. Suppose that $\om$ is
a faithful normal state on $N$ that is invariant under the inner
action
$$\be : N \recht M \ot N : \be(z) = U^* (1 \ot z) U \; ,$$
i.e., $(\io \ot \om)\be(z) = \om(z) 1$ for all $z \in N$.
Writing $X_n$ for the $n$-fold tensor product, $X_n:=U_{1,n+1} \ldots
U_{12} \in M \ot N^{\ot n}$, we define the inner action
$$\be_n : N^{\ot n} \recht M \ot N^{\ot n} : \be_n(z) = X_n^* (1 \ot
z) X_n \; .$$ It is clear that the tensor product state $\om^{\ot n}$
is invariant under $\be_n$ and that the actions $\be_n$ and
$\be_{n+1}$ are compatible with the inclusion $N^{\ot n}
\hookrightarrow N^{\ot (n+1)} : z \mapsto z \ot 1$. So, we can take
easily the direct limit $(N_\infty,\om_\infty) = \bigotimes_1^\infty
(N,\om)$ with the limit action $\al : N_\infty \recht M \ot N_\infty$
that we call the infinite tensor product action.

We will give a sufficient condition for $\al$ to be a strictly outer
action.

Recall that the l.c.\ quantum group $(M,\de)$ has a universal
C$^*$-algebraic dual $(\Ahu,\dehu)$, such that the
$^*$-representations of $\Ahu$ are in one-to-one correspondence with
the unitary corepresentations of $(M,\de)$. We denote by $\cW \in \M(A
\ot \Ahu)$ the universal corepresentation. Here $A$ is the reduced
C$^*$-algebra of $(M,\de)$, which can be defined as the norm closure
of $\{(\io \ot \om)(W) \mid \om \in \B(H)_* \}$.
See \cite{JK} for details.

\begin{proposition} \label{prop.infinitetensor}
Let $(M,\de)$ be a locally compact quantum group. Let $N$ be a factor
with a faithful state $\om$ and let $U \in M \ot N$ be a unitary
corepresentation of $(M,\de)$ in $N$. Suppose that the following two
conditions hold.
\begin{enumerate}
\item The state $\om$ is invariant under the inner action 
$\be : N \recht M \ot N : \be(z) = U^* (1 \ot z) U$.
\label{cond.one}
\item Let $\rho : \Ahu \recht N$ be the $^*$-representation satisfying
  $(\io \ot \rho)(\cW) = U$. Then, every bounded sequence $(a_n)$ in
  $\M(\Ahu)$ satisfying $(\io \ot \om\rho)\dehu(a_{n+1}) = a_n$ for
  all $n$, is a constant scalar sequence.
\label{cond.two}
\end{enumerate}
Then, the infinite tensor product action $\al : N_\infty \recht M \ot
N_\infty$ is strictly outer.
\end{proposition}

\begin{proof}
Denote by $E_n : N_\infty \recht N^{\ot n}$ the natural conditional
expectations. We have $(\io \ot E_n)\al(z) = \be_n(E_n(z))$ for all $z
\in N_\infty$.

Since $\be_n$ is an inner action of $(M,\de)$ on $N^{\ot n}$ (i.e.\ cocycle equivalent with the
trivial action), we obtain that $\Ad X_n : M \kruisje{\be_n} (N^{\ot n}) \recht \Mh \ot N^{\ot n}$
is an isomorphism, sending $\be_n(N^{\ot n})$ to $1 \ot N^{\ot n}$. So, we conclude that
$$M \kruisje{\be_n} (N^{\ot n}) \cap \be_n(N^{\ot n})' = X_n^* (\Mh
\ot 1) X_n \; .$$

Let $z \in M \kruisje{\al} N_\infty \cap \al(N_\infty)'$.
Considering $M \kruisje{\al} N_\infty$ as a subalgebra of $\B(H) \ot N_\infty$, we can write $z_n = (\io \ot E_n)(z)$.
Then, for all $n$, $z_n \in M \kruisje{\be_n} (N^{\ot n}) \cap \be_n(N^{\ot n})'$. So, we can take $y_n \in \Mh$, with $\|y_n\| \leq \|z\|$, such
that
$$(\io \ot E_n)(z) = X_n^* (y_n \ot 1) X_n \; .$$
Since $E_n \circ E_{n+1} = E_n$, we get
that $y_n = (\io \ot \om)(U^* (y_{n+1} \ot 1) U)$ for all $n$.

Define the (right) C$^*$-algebraic action $\eta: \cK(H) \recht
\M(\cK(H) \ot \Ahu) : \eta(y) = \cW^* (y \ot 1) \cW$. Fix $\mu \in
\B(H)_*$ and define $a_n = (\mu \ot \io)\eta(y_n)$. Clearly, $\|a_n\|
\leq \|\mu\| \; \|z\|$ for all $n$ and
$$(\io \ot \om\rho)\dehu(a_{n+1}) = (\mu \ot \io \ot \om\rho)(\io \ot
\dehu)\eta(y_{n+1}) = (\mu \ot \io)\eta\bigl( (\io \ot
\om\rho)\eta(y_{n+1}) \bigr) = (\mu \ot \io)\eta(y_n) = a_n \; .$$
By assumption, the sequence $(a_n)$ is a constant scalar
sequence. Since this holds for all $\mu \in \B(H)_*$, we conclude that
$\eta(y_n) \in \B(H) \ot 1$. But then, $\eta(y_n) = y_n \ot 1$. This
implies that $y_n \in M'$. We already know that $y_n \in \Mh$ and
conclude that $y_n \in \C 1$ for all $n$. Hence, $z \in \C 1$ and we
are done.
\end{proof}

In certain cases, the conditions in Proposition
\ref{prop.infinitetensor} can be weakened:
\begin{proposition} \label{prop.easierconditions}
Suppose that we are in the setting of Proposition \ref{prop.infinitetensor}.
\begin{itemize}
\item If the scaling group $(\tau_t)$ of $(M,\de)$ is trivial and
  $\om$ is a trace, then condition \ref{cond.one}) is automatically
  fulfilled.
\item If there exists a state $\om_1$ on $\Ahu$ and a number $0 < t <
  1$ such that $\om \rho = (1-t) \om_1 + t \epsh$, where $\epsh$ is
  the co-unit of $(\Ahu,\dehu)$, then condition \ref{cond.two}) can be
  weakened to the condition: every $a \in \M(\Ahu)$ satisfying $(\io
  \ot \om\rho)\dehu(a) = a$ is scalar.
\end{itemize}
\end{proposition}

\begin{proof}
First, suppose that the scaling group of $(M,\de)$ is trivial and that
$\om$ is a trace. Then,
also the scaling group of $(\Ahu,\dehu)$ is trivial and we get a
$^*$-anti-automorphism $\Rhu$ of $\Ahu$ such that
$$\Rhu \bigl( (\mu \ot \io)(\cW) \bigr) = (\mu \ot \io)(\cW^*)$$
for all $\mu \in \B(H)_*$. Let $f,g \in H$ and let $(e_n)$ be an
orthonormal basis for $H$. Let $z \in N$. We make the following
computation, writing $\om_{f,g} \in \B(H)_*$ defined by $\om_{f,g}(x)
= \langle x f , g \rangle$.
\begin{align*}
(\om_{f,g} \ot \om)(U^* (1 \ot z) U) &= \om \bigl( (\om_{f,g} \ot
\io)(U^* (1 \ot z) U) \bigr) = \sum_{n = 1}^\infty \om \bigl(
(\om_{e_n,g} \ot \io)(U^*) z  (\om_{f,e_n} \ot \io)(U) \bigr) \\ &=
\sum_{n=1}^\infty \om \bigl( z \; \rho\bigl( (\om_{f,e_n} \ot
\io)(\cW) (\om_{e_n,g} \ot \io)(\cW^*) \bigr) \bigr) \; .
\end{align*}
We now use that, with strict convergence,
\begin{align*}
\sum_{n=1}^\infty (\om_{f,e_n} \ot
\io)(\cW) (\om_{e_n,g} \ot \io)(\cW^*) &= \Rhu\Bigl( 
\sum_{n=1}^\infty (\om_{e_n,g} \ot \io)(\cW) (\om_{f,e_n} \ot
\io)(\cW^*) \Bigr) \\ &= \Rhu( (\om_{f,g} \ot \io)(\cW \cW^*)) =
\om_{f,g}(1) 1 \; .
\end{align*}
We conclude that
$$(\om_{f,g} \ot \om)(U^* (1 \ot z) U) = \om_{f,g}(1) \om(z)$$
which yields a proof of our first statement.

Suppose next that $\om \rho = (1-t) \om_1 + t \epsh$
Let $(a_n)$ be a bounded sequence in
  $\M(\Ahu)$ satisfying $(\io \ot \om\rho)\dehu(a_{n+1}) = a_n$ for
  all $n$.

Define the probability measure $\gamma$ on $\Z$ by
  $\gamma(\{0\}) = t$ and $\gamma(\{1\}) = 1-t$. If we define the
  $n$-fold convolution $\om_n = \om_1 * \cdots * \om_1$ ($n$ times),
  we get
$$(\om\rho)^{* n} = \sum_{k=0}^\infty \gamma^{* n}(\{k\}) \om_k \; ,$$
where $\gamma^{* n}$ is the $n$-fold convolution of the measure
$\gamma$ on $\Z$. By Corollary 2 in \cite{foguel}, we know that
$\|\gamma^{* n} - \gamma^{* (n+1)}\|_1 \recht 0$ if $n \recht
\infty$. So, we find that $\|(\om\rho)^{* n} - (\om\rho)^{* (n+1)}\|
\recht 0$ if $n \recht \infty$. Then, we have, for all $n$ and $k$,
\begin{align*}
\|a_{k+1} - a_{k} \| &= \| (\io \ot (\om\rho)^{* n})\dehu(a_{n+k+1}) -
(\io \ot (\om\rho)^{* (n+1)})\dehu(a_{n+k+1}) \| \\ &\leq \|(\om\rho)^{*
  n} - (\om\rho)^{* (n+1)}\| \; \|a_{n+k+1}\| \; .\end{align*}
We let $n \recht \infty$ and use that the sequence $(a_n)$ is bounded, to
conclude that $a_{k+1} = a_k$ for all $k$. Hence, there exists an $a
\in \M(\Ahu)$ such that $a_k=a$ for all $k$. Then, $(\io \ot
\om\rho)\dehu(a) = a$ and our weakened condition yields that $a \in \C
1$. Hence, $(a_n)$ is a constant scalar sequence.
\end{proof}

As a final general result, we prove that the possibility of acting strictly outerly on an injective factor is stable under cocycle perturbation. More
precisely, let $(M,\de)$ be a l.c.\ quantum group. A $2$-cocycle $\Om$ on $M$ is a unitary operator $\Om \in M \ot M$ satisfying
\begin{equation} \label{eq.cocycle}
(1 \ot \Om) (\io \ot \de)(\Om) = (\Om \ot 1) (\de \ot \io)(\Om) \; .
\end{equation}
Let as in the Preliminaries, $W$ denote the left regular corepresentation of $(M,\de)$, which is a multiplicative unitary on $H \ot H$. Recall that we
introduced as well the notations $J$ and $\Jh$ for the modular conjugations of the left invariant weights on $(M,\de)$ and $(\Mh,\deh)$,
respectively.

Suppose that $\Om$ is a $2$-cocycle on $(M,\de)$. Following \cite{Vai}, we define
$$\Wom := \Omtil W \Om^* \quad\text{with}\quad \Omtil = (1 \ot J \Jh) \Si \Om \Si (1 \ot \Jh J) \; .$$
We also define $\deom(x) = \Om \de(x) \Om^*$. It is easy to verify that $\Wom$ is a multiplicative unitary on $H \ot H$ and that $\deom$ is a
co-associative comultiplication on $(M,\de)$. Up to now, a general theory giving necessary and sufficient conditions for $(M,\deom)$ to be a l.c.\
quantum group is not available. See \cite{Vai} for results in this direction.
Nevertheless, we have the following result.

\begin{proposition} \label{prop.cocycle}
Let $(M,\de)$ be a l.c.\ quantum group and $\al : N \recht M \ot N$ a strictly outer action on the injective factor $N$. If $\Om$ is a $2$-cocycle on
$M$ such that $(M,\deom)$ is again a l.c.\ quantum group with left regular corepresentation $\Wom$, then
$$\be : N \ot \B(H) \recht M \ot N \ot \B(H) : \be(z) = (\Om W^*)_{13} \; (\al \ot \io)(z) \; (W \Om^*)_{13}$$
defines a strictly outer action of $(M,\deom)$ on the injective factor $N \ot \B(H)$.
\end{proposition}
\begin{proof}
It is obvious that $\al : N \recht M \ot N$ is a cocycle action (in the sense of \cite{VV}, Definition 1.1) of $(M,\deom)$ on $N$ with cocycle $\Om^*
\ot 1$. Using Definition 1.3 in \cite{VV}, we can define the cocycle crossed product
$$(M,\deom) \, \kruisje{(\al,\Om^* \ot 1)} N := \bigl( \; \al(N) \; \cup \; \{ (\om \ot \io)(\Omtil W) \ot 1 \mid \om \in \B(H)_* \} \; \bigr)\dpr \; .$$
If we consider the cocycle action $\al \ot \io$ of $(M,\deom)$ on $N \ot \B(H)$, with cocycle $\Om^* \ot 1 \ot 1$, the cocycle formula
\eqref{eq.cocycle} yields that this cocycle action is stabilizable with the unitary $(\Om W^*)_{13} \in M \ot N \ot \B(H)$, in the sense of
\cite{VV}, Definition 1.7. This means that, defining $\be$ as in the statement of the proposition, we indeed get an action of $(M,\deom)$ on $N \ot
\B(H)$. We claim that this action is strictly outer. Using Proposition 1.8 of \cite{VV}, we have to prove that
$$(M,\deom) \, \kruisje{(\al,\Om^* \ot 1)} N \; \cap \; \al(N)' = \C \; .$$
The left hand side is a subalgebra of $\B(H) \ot N \cap \al(N)' = M'
\ot 1$ by the strict outerness of $\al$ and Proposition \ref{prop.alternativeouter}.

As an intermediate step, we show that every $z \in (M,\deom) \,
\kruisje{(\al,\Om^* \ot 1)} N$ satisfies
\begin{equation}\label{eq.intermediate}
(\io \ot \al)(z) = \bigl( V \; (J \Jh \ot 1) \Si \Om^* \Si (\Jh J \ot 1) \bigr)_{12} \; z_{13} \; \bigl( V \; (J \Jh \ot 1) \Si \Om^* \Si (\Jh J \ot
1) \bigr)_{12}^* \; ,
\end{equation}
where $V \in \Mh' \ot M$ denotes the right regular corepresentation of $(M,\de)$. Recall that $$V = (J \Jh \ot 1) \Si W \Si (\Jh J \ot 1)$$ and $\de(x) =
V(x \ot 1)V^*$ for all $x \in M$.
Since, for $x \in N$, $(\io \ot \al)\al(x) = (\de \ot \io)\al(x) = V_{12} \al(x)_{13} V_{12}^*$ and since $(J \Jh \ot 1) \Si \Om^* \Si (\Jh J \ot 1)$
belongs to $M' \ot M$, Equation \eqref{eq.intermediate} is clear for $z = \al(x)$. To obtain Equation \eqref{eq.intermediate} for $z = (\om \ot
\io)(\Omtil W) \ot 1$ and $\om \in \B(H)_*$, we have to show that
$$\bigl( V \; (J \Jh \ot 1) \Si \Om^* \Si (\Jh J \ot 1) \bigr)_{23} \quad\text{and}\quad (\Omtil W)_{12} \quad\text{commute.}$$
Using the equality for $V$ recalled above, we have to show that $(\Om W^*)_{12}$ and $(\Om V)_{23}$ commute. Observe that
$$(\Om W^*)_{12} \; (\Om V)_{23} \; (W \Om^*)_{12} = \Om_{12} \; (\de \ot \io)(\Om) \; V_{23} \; \Om^*_{12} = \Om_{23} \; (\io \ot \de)(\Om) \; V_{23} \;
\Om^*_{12}= (\Om V)_{23} \; ,$$ where we used the cocycle equation \eqref{eq.cocycle}. Hence, we have proven Equation \eqref{eq.intermediate}.

Let now $z \in (M,\deom) \, \kruisje{(\al,\Om^* \ot 1)} N \; \cap \; \al(N)'$. Then, $z = J \Jh a \Jh J \ot 1$ for $a \in M$. Since $z$ satisfies Equation
\eqref{eq.intermediate}, we conclude that $a \ot 1$ and $\Si W \Om^* \Si$ commute. Hence, $\deom(a) = 1 \ot a$. So, $a \in \C$ and we are done.
\end{proof}

\section{The case of compact quantum groups}

As another partial converse to Proposition \ref{prop.amenable}, we
show that every compact Kac algebra whose dual has an invariant mean,
can act strictly outerly on the hyperfinite II$_1$ factor. In \cite{HY}, Theorem 8.6, the same kind of result is stated,
but we were unable to understand their complicated approach. Below we
present a fairly easy construction.

At the end of this section, we show that it is highly improbable that
there exists, for all $0 < q < 1$, a strictly outer action of the compact quantum group
$SU_q(2)$ on an injective factor.

Recall that we still continue to assume that our quantum groups are
second countable (i.e.\ all Hilbert spaces and preduals of von Neumann
algebras are separable).

In our von Neumann algebraic setting, a l.c.\ quantum group $(M,\de)$ is compact if and only if its Haar measure is finite. We refer to \cite{MVD}
for a nicely written overview of the theory of compact quantum groups
and their duals: discrete quantum groups. Suppose that $(M,\de)$ is a compact
quantum group. There exists a unique left invariant state $h$ on $M$, which is called the Haar state. The state $h$ is right invariant as well (a
compact quantum group is unimodular). Denote by $(\Mh,\deh)$ the dual l.c.\ quantum group, constructed as in the preliminaries out of the
multiplicative unitary $W \in M \ot \Mh$. We know that $\Mh = \bigoplus_{\al \in I} \Mh_\al$, where $\Mh_\al$ are finite-dimensional full matrix
algebras. Denote by $e_\al$ the minimal central projections of $\Mh$. Then, $U_\al := W(1 \ot e_\al)$ provide exactly the irreducible
corepresentations of $(M,\de)$. A compact quantum group is a compact
Kac algebra if and only if the scaling group $(\tau_t)$ is trivial. Equivalently, the
Woronowicz characters on the Hopf subalgebra $\cA \subset M$ (of matrix coefficients of finite dimensional corepresentations) are trivial.

The following lemma is a quantum version of Theorem 4.3 in \cite{KaiVer}. The result only works in the Kac case. Even a normal $q$-trace $\phi$ which
satisfies the same conclusion does not exist in the non-Kac case.

Recall that we define the convolution product $\om * \mu := (\om \ot \mu)\deh$ for $\om,\mu \in \Mh_*$.

\begin{lemma} \label{lem.tracestate}
Let $(\Mh,\deh)$ be a discrete Kac algebra with an invariant
mean. Then, there exists a normal tracial
state $\phi$ on $\Mh$ such that
$$\| \mu * \phi_n - \mu(1) \, \phi_n \| \recht 0 \quad\text{for
  all}\quad \mu \in \Mh_* \quad\text{where}\quad \phi_n =
\underset{n \; \text{times}}{\underbrace{\phi * \cdots * \phi}} \; .$$
\end{lemma}
\begin{proof}
As a first step, we prove that there exists a sequence $\om_n$ of normal tracial
states on $\Mh$ such that $\| \mu * \om_n - \mu(1) \, \om_n \| \recht
0$ for all $\mu \in \Mh_*$. Let $m_0$ be an invariant mean on
$(\Mh,\deh)$. Denote by $W \in M \ot \Mh$ the multiplicative unitary
of $(M,\de)$. Define
$$\Phi : \Mh \recht M \ot \Mh : \Phi(x) = W^*(1 \ot x)W \; .$$
Then, $\Phi$ is an action of $(M,\de)$ on $\Mh$ (the adjoint
action). Let $h$ be the Haar state on the compact quantum group
$(M,\de)$ and define a (non-normal) state $m$ on $\Mh$ by the formula
$$m(x) = m_0\bigl( (h \ot \io)\Phi(x) \bigr) \quad\text{for}\quad x
\in \Mh \; .$$
We claim that $m$ is a left invariant mean on $(\Mh,\deh)$. Let $x \in
\Mh$ and $\mu \in \Mh_*$. Then,
\begin{align*}
m \bigl( (\io \ot \mu)\dehop(x) \bigr) &= m_0\bigl( (h \ot \io \ot
\mu)(W^*_{12} W_{23} (1 \ot x \ot 1) W_{23}^* W_{12} ) \bigr) \\ &=
m_0\bigl( (h \ot \io \ot \mu)(W_{13} W_{23} (\Phi(x) \ot 1) W_{23}^*
W_{13}^* ) \bigr) \\ &= m_0\bigl( (((h \ot \mu) \circ \Ad W) \ot \io)(\io \ot
\deh)\Phi(x)\bigr) \\ &= m_0\bigl( (((h \ot \mu) \circ \Ad W) \ot
\io)(\Phi(x)_{13}) \bigr) \; .
\end{align*}
It is easy to check that, for all $z \in M$, $(h \ot \io)(W(z \ot 1)W^*) = h(z) \, 1$, see e.g.\ the proof of Corollary 3.9 in \cite{Iz}. Here, we
have used in a crucial way that we are in the Kac case. So, we can continue the computation above and get
$$m \bigl( (\io \ot \mu)\dehop(x) \bigr) = \mu(1) \; m_0\bigl( (h \ot
\io)\Phi(x) \bigr) = \mu(1) \; m(x) \; .$$
Hence, $m$ is indeed a left invariant mean on $(\Mh,\deh)$.

Take a sequence $\eta_n$ of normal states on $\Mh$ such that
$\eta_n(x) \recht m_0(x)$ for all $x \in \Mh$. Define $\ga_n := (h \ot
\eta_n) \Phi$. Then, $\ga_n$ is a sequence of normal states on $\Mh$
such that $\ga_n(x) \recht m(x)$ for all $x \in \Mh$. By invariance of
$h$, we also have that, for all $n$, $\ga_n$ is invariant under the
action $\Phi$: $(\io \ot \ga_n)\Phi(x) = \ga_n(x) \; 1$ for all $x \in
\Mh$.

We claim that a normal state $\om$ on $\Mh$ is a trace if and only if $\om$ is
invariant under $\Phi$. We know that $\Mh$ is a direct sum of matrix
algebras $\Mh_\al$. It is also clear that $\Phi$ leaves invariant all
the matrix algebras $\Mh_\al$, yielding actions $\Phi_\al$ of
$(M,\de)$ on $\Mh_\al$.
Every of these restrictions $\Phi_\al$ is
ergodic (i.e.\ has a trivial fixed point algebra).
Let $\tau_\al$ be the normalized trace on
$\Mh_\al$. It is easy to verify that $\tau_\al$ is invariant under
$\Phi_\al$ (see e.g.\ Lemma 2.1 in \cite{Iz}). If now $\om$ is a
normal state on $\Mh$ which is invariant under $\Phi$, we get that
$\om = \sum \tau_\al(K_\al \cdot)$, for certain positive matrices
$K_\al$, satisfying $\Phi_\al(K_\al) = 1 \ot K_\al$. Hence, all the
$K_\al$ are scalar and $\om$ is a trace.

We conclude that $\ga_n$ is a sequence of normal tracial states on
$\Mh$. Since $\ga_n \recht m$ pointwise and since $m$ is left
invariant, a classical technique allows to find a sequence of normal
states $\om_n$ that are convex combinations of the $\ga_n$ and that
satisfy $\| \mu * \om_n - \mu(1) \, \om_n \| \recht 0$ for all $\mu
\in \Mh_*$. Since a convex combination of traces is still a trace, the
first step of the proof is done.

For the second step of the proof, we can follow almost literally the
proof of Theorem 4.3 in \cite{KaiVer}, yielding a sequence $t_n$ of
positive real numbers, satisfying $\sum t_n = 1$ and such that
$$\phi:= \sum_n t_n \om_n$$
is the normal tracial state that we are looking for.
\end{proof}

This allows us to prove the announced result.

\begin{theorem} \label{thm.compactKac}
Let $(M,\de)$ be a compact Kac algebra. If the dual discrete Kac
algebra $(\Mh,\deh)$ has an invariant mean, then $(M,\de)$ can act
strictly outerly on the hyperfinite II$_1$ factor.
\end{theorem}
\begin{proof}
Take a normal tracial state $\phi$ satisfying the conclusion of Lemma \ref{lem.tracestate}. Write $\phi = \sum t_\al \tau_\al$, where $\tau_\al$ is
the normalized trace on $\Mh_\al$. Take the subset $I_0$ of $\al \in
I$ satisfying $t_\al > 0$. Denote by $\cR$ the hyperfinite II$_1$ factor and
denote by $\tau$ its tracial state. Take a family of orthogonal
projections $\{ e_\al  \in \cR \mid \al \in I_0 \}$ such that $\tau(e_\al) = t_\al$.
Choosing isomorphisms $e_\al \cR e_\al \cong \Mh_\al \ot \cR$, we get embeddings $\Mh_\al \rightarrow e_\al \cR e_\al$ and hence, a normal
$^*$-homomorphism $\pi : \Mh \recht \cR$. Define $U := (\io \ot
\pi)(W) \in M \ot \cR$ and $\be : \cR \recht M \ot \cR : \be(x) = U^*
(1 \ot x) U$. From Proposition \ref{prop.easierconditions}, we know
that the trace $\tau$ is invariant under $\be$. So, as explained
before Proposition \ref{prop.infinitetensor}, we can make an
infinite tensor product action $\al$ of $(M,\de)$ on $\cR \cong
\bigotimes_1^\infty \cR$. In order to prove that $\al$ is strictly
outer, we have to verify condition \ref{cond.two}) in Proposition
\ref{prop.infinitetensor}. Observe that, since $(M,\de)$ is compact,
$\M(\Ahu) = \Mh$. From Lemma \ref{lem.tracestate}, we know
that $\|\phi_n - \phi_{n+1}\| \recht 0$. Hence, as in the proof of the
second item of Proposition \ref{prop.easierconditions}, it suffices to
check that an element $a \in \Mh$ satisfying $(\io \ot \phi)\deh(a) =
a$ is scalar. But, such an element satisfies $\mu(a) = (\mu *
\phi_n)(a)$ for all $n$ and hence, $|\mu(a) - \mu(1) \phi_n(a)| \recht
0$ if $n \recht \infty$. This means that $a - \phi_n(a) 1 \recht 0$
weakly and hence, $a$ is scalar.
\end{proof}

We now prove, combining results of Ueda \cite{U2} and a classical trick, that a strictly outer action of $SU_q(2)$ on an injective factor gives rise to
an irreducible subfactor of the hyperfinite II$_1$ factor with index
$(q+q^{-1})^2$.

We will make use of the extension of Jones' index theory to factors which are not necessarily type II$_1$, see \cite{kos}. In this theory, an index
$\ind E$ is associated to an inclusion $N \subset M$ and a faitfhul normal conditional expectation $E : M \recht N$. It coincides with the Jones
index, if $N$ and $M$ are II$_1$ factors and $E$ is the unique conditional expectation satisfying $\tau_N \circ E = \tau_M$ (where $\tau_N$ and
$\tau_M$ are the normalized traces on $N$ and $M$, resp.).

The following result is probably well known, but we include a proof
for the convenience of the reader. It uses essentially a non-compact version of Wassermann's
invariance principle \cite{Was}.

\begin{lemma} \label{lem.reduce}
Let $N \subset M$ be an inclusion of factors and $E : M \recht N$ a faithful normal conditional expectation with $\ind E = \sla < \infty$. Let $N
\subset M \subset M_1 \subset \cdots$ be the Jones tower with associated conditional expectations $E_n : M_n \recht M_{n-1}$.

Then, there exists an inclusion $\Ntil \subset \Mtil$ of II$_1$
factors with $[\Mtil : \Ntil] = \sla$ and whose Jones tower $\Ntil \subset \Mtil \subset
\Mtil_1 \subset ...$ satisfies
$$\Mtil_i \cap \Mtil_{j-1}' = (M_i \cap M_{j-1}')^{(\si_t^{E_{ij}})} \quad\text{where}\quad E_{ij} = E_j \circ E_{j+1} \circ \cdots \circ E_i \; .$$

If $M$ is injective, $\Ntil \cong \Mtil \cong \cR$, the hyperfinite II$_1$ factor.
\end{lemma}
\begin{proof}
Take a faithful normal state $\eta$ on $N$. Put $\eta_0 := \eta E$ and $\eta_i := \eta_{i-1} E_i$. Let $A$ be the injective factor of type III$_1$,
with faithful normal state $\mu$. Write $\om_i = \eta_i \ot \mu$ on the III$_1$ factor $M_i \ot A$ and $\om = \eta \ot \mu$ on $N \ot A$. Consider
the cores $C_i := (M_i \ot A) \rtimes_{(\si_t^{\om_i})} \R$ and $C:= (N \ot A) \rtimes_{(\si_t^{\om})} \R$, which we realize as subalgebras of $M_i
\ot A \ot \B(L^2(\R))$, respectively. Since $M_i \ot A$ is a III$_1$
factor, $C_i$ is a II$_\infty$ factor for all $i$.

The restriction of $E_i \ot \io \ot \io$ yields a conditional expectation $F_i : C_i \recht C_{i-1}$ and $F : C_0 \recht C$. By the characterization
of the Jones tower (\cite{HK}, Theorem 8), $C \subset C_0 \subset C_1 \subset \cdots$ is a Jones tower with compatible conditional expectations
$F,F_i$. In particular, $\ind F = \sla$. Denote by $\omtil$ and $\omtil_i$ the dual weights on $C$ and $C_i$, resp. It is clear that $\omtil_i =
\omtil_{i-1} F_i$. If we denote by $\rho_t$ the right translation operators in $\B(L^2(\R))$, we can consider $1 \ot 1 \ot \rho_t \in C_i$. We can
define traces $\Trace_i$ on the II$_\infty$ factors $C_i$ (and $\Trace$ on $C$) such that the Connes cocycles w.r.t.\ $\omtil_i$ are given by
$$[\omtil_i : \Trace_i]_t = 1 \ot 1 \ot \rho_t \; .$$
It follows that $\Trace_i = \Trace_{i-1} F_i$. Let $p$ be a projection in $C$ satisfying $\Trace(p) = 1$. Put $\Ntil = pCp$, $\Mtil = pC_0 p$ and
$\Mtil_i = pC_i p$. It is clear that we get conditional expectations $\Etil_i : \Mtil_i \recht \Mtil_{i-1}$ and $\Etil : \Mtil \recht \Ntil$ such
that the towers
\begin{equation*}
C \overset{\overset{F}{\leftarrow}}{\subset} C_0 \overset{\overset{F_1}{\leftarrow}}{\subset} C_1 \overset{\overset{F_2}{\leftarrow}}{\subset} \cdots
\quad\text{and}\quad \Ntil \ot \B(\ell^2(\Z)) \overset{\overset{\Etil \ot \io}{\leftarrow}}{\subset} \Mtil \ot \B(\ell^2(\Z))
\overset{\overset{\Etil_1 \ot \io}{\leftarrow}}{\subset} \Mtil_1 \ot \B(\ell^2(\Z)) \overset{\overset{\Etil_2 \ot \io}{\leftarrow}}{\subset} \cdots
\end{equation*}
are isomorphic in a way that preserves the conditional expectations. Since $\Tr_i(p) = 1$, it follows that $\Ntil \subset \Mtil$ is an inclusion of
II$_1$ factors. It is clear that $\ind \Etil = \sla$. Because $\Tr
\circ F = \Tr_0$ and because the restriction of $\Tr$ to $\Ntil$ and of $\Tr_0$ to
$\Mtil$ are the unique tracial states of $\Ntil$ and $\Mtil$, we get $[\Mtil : \Ntil]=\sla$.

Next, we observe that
$$\Mtil_i \cap \Mtil_{j-1}' \cong (M_i \ot A) \rtimes_{(\si_t^{\om_i})} \R \cap ((M_{j-1} \ot A) \rtimes_{(\si_t^{\om_i})} \R)' =C_i \cap C_{j-1}' \; .$$
We compute
$$C_i \cap C_{j-1}' \subset M_i \ot A \ot \B(L^2(\R)) \cap (1 \ot
\si^\mu(A))' = M_i \ot 1 \ot L^\infty(\R)$$ because $A$ is a III$_1$ factor. Hence, we get
$$C_i \cap C_{j-1}' \subset \si^{\eta_i}(M_i)_{13} \cap
\si^{\eta_{j-1}}(M_{j-1})_{13}' \cap (1 \ot 1 \ot \rho(\R))' = (M_i \cap M_{j-1}')^{(\si_t^{E_{ij}})} \ot 1 \ot 1 \; .$$ Since the converse inclusion
$(M_i \cap M_{j-1}')^{(\si_t^{E_{ij}})} \ot 1 \ot 1 \subset C_i \cap C_{j-1}'$ is clear, we have proven the formula for $\Mtil_i \cap \Mtil_{j-1}'$.

To conclude the proof, we see that if $M$ is injective, $\Mtil$ and $\Ntil$ are injective factors of type II$_1$. Hence, $\Ntil \cong \Mtil \cong
\cR$ in that case.
\end{proof}

Combining the previous lemma with results of Ueda \cite{U2}, we get the following proposition. Let $(M,\de)$ be a compact quantum group and $u \in M
\ot \M_n(\C)$ an irreducible corepresentation on $\C^n$. Associated
with such an irreducible corepresentation is a positive invertible $F$-matrix $F_u
\in \M_n(\C)$ satisfying $\Trace F_u = \Trace F_u^{-1}$ (see \cite{Wor1}). The positive real number $\Trace F_u$ is denoted by $\dim_q u$ and called the quantum
dimension of $u$.

\begin{proposition} \label{prop.irreducible}
Let $\al : N \recht M \ot N$ be a strictly outer action of a compact quantum group $(M,\de)$ on an injective factor $N$. Then, there exist
irreducible subfactors of the hyperfinite II$_1$ factor with index $(\dim_q u)^2$ for any irreducible corepresentation $u$ of $(M,\de)$.
\end{proposition}
\begin{proof}
Let $\al$ be such an action and $u \in M \ot \M_n(\C)$ an irreducible corepresentation. Recall that for compact quantum groups the notions of
strictly outer action and minimal action coincide. We can define a new action
$$\ga : N \ot \M_n(\C) \recht M \ot N \ot \M_n(\C) : \ga(z) = u^*_{13} \; (\al \ot \io)(z) \; u_{13} \; .$$
Consider the inclusion $N^\al \ot 1 \subset \bigl( N \ot \M_n(\C)\bigr)^\ga$. Using the restriction $E_u$ of $\frac{1}{\dim_q u}(\io \ot \Tr(F_u
\cdot))$, Ueda proved in \cite{U2} that $N^\al \ot 1 \subset \bigl( N
\ot \M_n(\C)\bigr)^\ga$ is an irreducible inclusion of factors with $\ind E_u
= (\dim_q u)^2$. Using Lemma \ref{lem.reduce}, we get the existence of an irreducible subfactor of the hyperfinite II$_1$ factor with index $(\dim_q
u)^2$.
\end{proof}

\begin{remark} \label{remark.suqtwo}
So, if the compact quantum group $SU_q(2)$ would have a strictly outer
action on an injective factor for all $0 < q < 1$, we can use its
fundamental corepresentation $u$ with $\dim_q u = q + q^{-1}$ and
conclude that the hyperfinite II$_1$ factor would have irreducible
subfactors of arbitrary index strictly greater than 4.
Since most important unpublished work of Popa states
that not all values strictly greater than $4$
can be realized as the index of an irreducible subfactor of the hyperfinite II$_1$ factor, it follows that at least for certain values of $0 < q <
1$, $SU_q(2)$ cannot act strictly outerly on an injective
factor. Remark that nevertheless, Banica showed \cite{ban2} that the
dual of $SU_q(2)$ has an invariant mean. This would mean that the
converse of Proposition \ref{prop.amenable} does not hold in its full generality.
\end{remark}

\begin{remark}
As we explained in the previous remark, there is a strong reason to believe that $SU_q(2)$ cannot act strictly outerly on an injective factor for certain values of
$q$. Recently, Szymanski \cite{szy} has shown that the compact quantum groups $SU_q(2)$ for different values of $q$ are related by a
pseudo-$2$-cocycle. If the pseudo-$2$-cocycle of Szymanski happens to
be a $2$-cocycle, we can use Proposition \ref{prop.cocycle} and can conclude that none of the compact quantum groups $SU_q(2)$, $0 < q < 1$, can act
strictly outerly on an injective factor.
\end{remark}

\begin{remark}
In the proof of Proposition \ref{prop.irreducible}, we constructed,
following Ueda \cite{U2}, an irreducible inclusion of factors given a
strictly outerly acting compact quantum group $(M,\de)$ and an irreducible
corepresentation. Ueda computed the Jones tower of this inclusion. The
tower of relative commutants only depends on the corepresentation
theory of the compact quantum group. The vertices of the principal
graph are labeled by the irreducible corepresentations of $(M,\de)$
that are subrepresentations of some tensor product $u \ot \overline{u}
\ot \cdots$. In particular, whenever the matrix coefficients of $u$
generate a C$^*$-subalgebra of $M$ of infinite dimension, we get an
inclusion of infinite depth. This is, of course, almost always the
case.
\end{remark}

\section{The case of discrete quantum groups}

Since the dual of a discrete quantum group is compact and since a
compact quantum group has an invariant mean (the Haar state),
Proposition \ref{prop.amenable} does not exclude the possibility that
every discrete quantum group can act strictly outerly on an injective
factor.

However, if a discrete quantum group $(\Mh,\deh)$ with invariant mean acts strictly
outerly on an injective factor, the crossed product will be injective
as well and the dual action will be a strictly outer action of the
compact quantum group $(M,\de)$. So, from Remark \ref{remark.suqtwo},
we conclude that we should not expect to find a strictly outer action
of the dual of $SU_q(2)$ on an injective factor.

We first prove that a discrete Kac algebra with invariant mean acts
strictly outerly on the hyperfinite II$_1$ factor $\cR$. Taking the crossed
product and the dual action, we get an alternative proof for Theorem
\ref{thm.compactKac}.

Next, we prove more generally that every discrete Kac algebra with a
faithful corepresentation in $\cR$ (Definition \ref{def.faithfulcorep}), acts strictly outerly on
$\cR$. This generalizes Banica's result (\cite{Ban}, Section 4) for
discrete Kac algebras with a faithful finite-dimensional
corepresentation. Note that the discrete quantum groups constructed
from vertex models have, by construction, such a faithful
finite-dimensional corepresentation.

Recall once again that we assume that our l.c.\ quantum groups are
second countable.

The following is an easy application of Propositions
\ref{prop.infinitetensor} and \ref{prop.easierconditions}.

\begin{proposition}
Let $(\Mh,\deh)$ be a discrete Kac algebra with invariant mean. Then,
there exists a strictly outer action of $(\Mh,\deh)$ on the
hyperfinite II$_1$ factor.
\end{proposition}

\begin{proof}
Let $h$ denote the Haar state on the compact Kac algebra
$(M,\de)$. Then, $h$ is a faithful tracial state. Moreover, since
$(\Mh,\deh)$ has an invariant mean, it is easy to check that $M$ is an
injective von Neumann algebra (see e.g.\ \cite{ruan}). Hence, there
exists a faithful $^*$-homomorphism $\rho : M \recht \cR$ such that
$\tau \rho = h$, where $\tau$ is the tracial state on $\cR$ (this
essentially follows from the uniqueness of $\cR$, see Corollary 1.23
in \cite{tak3}).

Consider the multiplicative unitary $\Wh = \Si W^* \Si \in \Mh \ot M$,
where $\Si$ denotes the flip map. Define $U:= (\io \ot \rho)(\Wh)$. We
claim that the conditions of Proposition \ref{prop.infinitetensor} are
fulfilled. Condition \ref{cond.one}) is fulfilled because of
Proposition \ref{prop.easierconditions}. On the other hand, because
$(\Mh,\deh)$ has an invariant mean, we know that $A\uni = A$ and
hence, $\M(A\uni) \hookrightarrow M$. If $(a_n)$ is a bounded sequence
in $M$ satisfying $(\io \ot \tau\rho)\de(a_{n+1}) = a_n$, we get that
$a_n = h(a_{n+1}) 1$, because $\tau \rho=h$. So, every $a_n$ is scalar
and then, $(a_n)$ is a constant scalar sequence.
\end{proof}

The non-trivial point in the proof of the previous proposition is the
existence of a faithful, normal $^*$-homomorphism $\rho : M \recht
\cR$. The existence of $\rho$ implies the injectivity of $M$ and this,
in turn, implies the existence of an invariant mean on $(\Mh,\deh)$
(see e.g.\ \cite{ruan}). So, we cannot follow the same strategy in the
non-amenable case.

We shall consider discrete Kac algebras that have a faithful
corepresentation in the hyperfinite II$_1$ factor $\cR$. In the classical
case, this corresponds to discrete groups that are subgroups of the
unitary group of $\cR$. In particular, all residually finite groups
belong to this class and they may very well be non-amenable.

\begin{theorem}
If a discrete Kac algebra $(\Mh,\deh)$ has a faithful corepresentation
$U \in \Mh \ot \cR$ in the hyperfinite II$_1$ factor $\cR$, then
$(\Mh,\deh)$ acts strictly outerly on $\cR$.
\end{theorem}
\begin{proof}
Take the $^*$-homomorphism $\rho : \Au \recht \cR$ such that $(\io \ot
\rho)(\cWh) = U$, where $\cWh \in \M(\Ah \ot \Au)$ is the universal
corepresentation of $(\Mh,\deh)$.
Taking the direct sum with the trivial corepresentation, we may assume that
there exists a projection $e_0 \in \cR \cap \rho(\Au)'$ satisfying $0
< \tau(e_0) < 1$ and $\rho(a) e_0 = \eps(a) e_0$ for all $a \in \Au$,
where $\eps$ denotes the co-unit of $(\Au,\deu)$. From the
faithfulness of $U$, we know that
$$^*\text{-alg} \; \{ (\io \ot \mu)(U) \mid \mu \in \cR_* \}
\quad\text{is weakly dense in} \;\; \Mh \; .$$
Taking the
tensor product of $U$ and its adjoint corepresentation, we may assume
that
\begin{equation}\label{eq.dichtheid}
\text{alg} \; \{ (\io \ot \mu)(U) \mid \mu \in \cR_* \}\quad\text{is
  weakly dense in} \;\; \Mh \; .
\end{equation}

From Proposition \ref{prop.easierconditions}, we know that $\tau$ is
invariant under the inner action $\be : \cR \recht \Mh \ot \cR :
\be(z) = U^*(1 \ot z)U$. So, we can define the infinite tensor product
action $\al$. We claim that $\al$ is strictly outer.

Write $\om = \tau\rho$. If we define the state $\om_1$ on $\Au$ such
that $(1 - \tau(e_0)) \om_1(a) = \tau(\rho(a)(1-e_0))$ and if we put
$t = \tau(e_0)$, we observe that $\om = (1-t)\om_1 + t \eps$. We apply
Propositions \ref{prop.infinitetensor} and
\ref{prop.easierconditions}. So, in
order to prove our claim, it suffices to check that the equation
$$(\io \ot \om)\de(a) = a \; , \quad a \in \M(\Au)$$
has only scalar solutions.

Let $\Psi$ be a state on $\Au$ that is an accumulation point in the
weak$^*$ topology of the sequence $(\Psi_n)$, where
$$\Psi_n := \frac{1}{n} \sum_{k=1}^n \om^{* k} \; .$$ We will
prove that $\Psi$ is the Haar state of $(\Au,\deu)$.
It is clear that $\om * \Psi = \Psi$. Define, for all $x \in \cR$, an
element $\mu_x \in \Au^*$ as
$\mu_x(a) = \tau(x^* \rho(a) x)$. Because $\mu_x \leq \|x\|^2 \om$, it
follows from Lemma 4.3 in \cite{MVD} that $\mu_x * \Psi = \mu_x(1)
\Psi$. So, we conclude that $(\mu\rho) * \Psi = \mu(1) \Psi$ for all
$\mu \in \cR_*$. Applying this to the second leg of $\cWh$, we get
$$(\io \ot \Psi)(\cWh) (\io \ot \mu\rho)(\cWh) = \mu(1) (\io \ot
\Psi)(\cWh)$$
for all $\mu \in \cR_*$. Denote by $\epsh \in \Mh_*$ the co-unit of
the discrete Kac algebra $(\Mh,\deh)$. Using Equation
\eqref{eq.dichtheid}, we conclude that $(\io \ot \Psi)(\cWh) a =
\epsh(a) (\io \ot \Psi)(\cWh)$ for all $a \in \Mh$. From this it
follows that $(\io \ot \Psi)(\cWh)$ is the central projection in $\Mh$
projecting on the trivial representation. Hence, $\Psi$ is the Haar
state.

If $a \in \M(\Au)$ and $(\io \ot \om)\de(a) = a$, then $(\io \ot
\Psi_n)\de(a) = a$ for all $n$. This implies that $(\io \ot
\Psi)\de(a) = a$. Since $\Psi$ is the Haar state of $(\Au,\deu)$, we
get $a = \Psi(a) 1$ and we are done.
\end{proof}

\end{document}